\documentclass[reqno]{amsart}
\usepackage{amscd}

\newtheorem{thm}{Theorem}[section]

\newtheorem{dfn}[thm]{Definition}

\newtheorem{prp}[thm]{Proposition}

\newtheorem{lma}[thm]{Lemma}
 
\newtheorem{rem}[thm]{Remark}
\newenvironment{rmk}{\begin{rem}\rm}{\end{rem}}

\newtheorem*{clm}{Claim}

\newenvironment{pf}{\begin{proof}}{\end{proof}}

\numberwithin{equation}{section}

\newcommand{\R}{{\mathbb{R}}}
\newcommand{\C}{{\mathbb{C}}}

\newcommand{\Z}{{\mathbb{Z}}}
\newcommand{\D}{{\mathbb{D}}}
\newcommand{\cle}{{\mathbb{S}}}
\newcommand{\cf}{{\sc cf}}
\newcommand{\cn}{{\sc cn}}
\newcommand{\rn}{{\sc rn}}
\newcommand{\rf}{{\sc rf}}
\newcommand{\rdp}{{\sc rdp}}
\newcommand{\sdp}{{\sc sdp}}
\newcommand{\po}{P}
\newcommand{\qu}{Q}
\newcommand{\la}{\langle}
\newcommand{\ra}{\rangle}
\newcommand{\pa}{\partial}
\newcommand{\lk}{\operatorname{lk}}
\newcommand{\wra}{\operatorname{wr}}
\newcommand{\sh}{\operatorname{sh}}
\newcommand{\EW}{\operatorname{{\mathbb C}{\mathit w}}}
\newcommand{\Hom}{\operatorname{Hom}}

\newcommand{\pr}{\operatorname{pr}}

\begin{document}


\title
{The complex shade of a real space and its applications}
\author{Tobias Ekholm}
\address{Department of Mathematics, Uppsala University, Box 480, 751
06 Uppsala, Sweden}
\email{tobias@math.uu.se}
\date{}






\thanks{During the preparation of this manuscript the author received
a post-doctoral grant from Stiftelsen f{\"o}r internationalisering av
forskning och h{\"o}gre utbildning.}


\begin{abstract}
A natural oriented $(2k+2)$-chain in $\C P^{2k+1}$
with boundary twice $\R P^{2k+1}$, its complex shade, is
constructed. Via intersection 
numbers with the shade, a new invariant, the shade number, of
$k$-dimensional subvarieties with normal vector fields along their
real part, is introduced. 
For an even-dimensional real variety, the shade number  
and the Euler number of the complement of the normal vector field in
the real normal bundle of its real part agree.
For an odd-dimensional orientable real variety, a linear
combination of the shade number and the wrapping number (self-linking
number) of its real part is independent
of the normal vector field and equals the encomplexed writhe as
defined by Viro ~\cite{V}. Shade 
numbers of varieties without real points and encomplexed 
writhes of odd-dimensional real varieties are, in a sense, Vassiliev  
invariants of degree $1$.

Complex shades of odd-dimensional spheres are constructed.
Shade numbers of real subvarieties in spheres have properties
analogous to those of their projective counterparts.   
\end{abstract}

\maketitle


\section{Introduction}\label{intr}
The manifold $\R P^n$, $n>0$, is orientable if and only if $n$ is odd.
With its standard orientation, $\R P^{2k+1}\subset\C P^{2k+1}$ is
homologous to zero, since $H_{2k+1}(\C P^{2k+1};\Z)=0$.
(Notice the difference between even and odd dimensions: 
$\R P^{2k}$ represents a generator of $H_{2k}(\C P^{2k};\Z_2)$. Our
main concern are subvarieties in ``double dimension $+$
$1$'' and so we consider only the odd-dimensional case.)
The {\em complex shade} of $\R P^{2k+1}$    
consists of the complexifications of all real lines through a
point $p$ in $\R P^{2k+1}$. It is an oriented $(2k+2)$-chain
$\Gamma_p$ in $\C P^{2k+1}$ with boundary twice $\R P^{2k+1}$, see
Definition ~\ref{dfnshadeP}. Geometrically distinct shades are in 1-1
correspondence with points in $\R P^{2k+1}$. They all represent the
same homology class, the {\em shade class} 
$[\Gamma]\in H_{2k}(\C P^{2k+1},\R P^{2k+1})$. 
This homology class $[\Gamma]$ has three characterizing properties:
it is invariant under complex conjugation, 
it has boundary twice $\R P^{2k+1}$, and  
its intersection number with the class represented by
a complex $k$-dimensional linear space without real points has
absolute value equal to $1$, see Proposition ~\ref{prpHax}. The shade class
will be used to measure certain linking phenomena.

Recall that if $M$ is an $m$-dimensional oriented manifold, and 
$B$ and $C$ are disjoint oriented cycles of dimensions $k$ and $m-k-1$,
respectively, which are weakly zero-homologous (i.e. represent
torsion classes in homology) then the {\em linking number} $\lk(B,C)$ is
defined as $\frac{1}{n}$ times the intersection number of an oriented 
$(k+1)$-chain $A$ with boundary $n\cdot B$ and $C$ in $M$.       

An oriented  $2k$-cycle $C$ in the complement of $\R P^{2k+1}$ in 
$\C P^{2k+1}$ need not be weakly zero-homologous and the linking
number of $\R P^{2k+1}$ and $C$ is in general not defined. 
In fact, there are many possible choices of a counterpart of linking
number in this case: the intersection number of any oriented
$2k$-chain with boundary $\R P^{2k+1}$ and $C$ would do. 
The shade class provides a choice: we measure the linking of $C$ and
$\R P^{2k+1}$ in $\C P^{2k+1}$ as $\frac12$ times the intersection
number of the shade class and $C$. This construction has
straightforward applications to $k$-dimensional complex projective
varieties without real points and give rise to what will be called shade
numbers, see \S ~\ref{IRempty}.  

Most $k$-dimensional complex varieties in projective
$(2k+1)$-space (varieties in open dense subset of the Chow-variety)
are without real points. But the definition of shade number can be
extended  to arbitrary $k$-dimensional projective subvarieties provided 
they have been equipped with additional structure. (The additional
structure is a certain kind of vector field along the real part in 
$\R P^{2k+1}$, on which the shade number will depend. Any
$k$-dimensional variety admits such additional structure and in
special cases it appears naturally, see \S ~\ref{secgeneral})   

The most interesting applications of shade numbers arise in the study
of the interplay between the real and the complex geometry/topology of 
generic real projective varieties (see \S ~\ref{Ishadearm}):  
there are connections between shade numbers 
and topological invariants of the real part
equipped with the above mentioned vector field. The nature of these
connections depends on the parity of the dimension of the variety and
the differences between even and odd dimensions are striking,
see \S ~\ref{Ieven} and \S ~\ref{Iodd}. Since
spheres frequently occurs as real algebraic ambient spaces, we 
define complex shades of spheres as well, see \S ~\ref{IRsphere}.

\subsection{Applications to real algebraic links}\label{applral}
{\em The writhe} of a generic
projection of a knot in $\R^3$ to a plane (a knot
diagram) is the signed sum of double points of its image. This notion
has a straightforward generalization to plane projections of
knots in $\R P^3$. In general the writhe changes as the projection varies.

It was observed by Viro ~\cite{V1} that the writhe of a knot may be
enhanced if it is real algebraic (a real algebraic knot in projective
$3$-space is a 
smooth 1-dimensional projective variety with connected real
part). In  ~\cite{V}, Viro introduced 
{\em the encomplexed writhe} which shows that the classification of real 
algebraic knots up to rigid isotopy (a smooth isotopy which is
also a continuous family of real algebraic curves, see ~\cite{Ro} and
~\cite{V2}) is more refined than the corresponding classification up
to smooth isotopy.   

Viro's definition is diagrammatic: the encomplexed writhe of a real
algebraic knot $V$
is expressed as the signed sum of real double points in a
generic real projection of $V$ to a plane. 
Note that the preimage of a double point in such a projection may be
either two distinct points in the real part of $V$,
or two complex conjugate non-real points of $V$. Double points of the
later kind are called solitary and they 
enhance the writhe to the encomplexed writhe, which is
independent of generic projection and invariant under rigid isotopy.    

The initial motivation for the study undertaken in this paper was to
give an intrinsic $3$-dimensional explanation of Viro's
invariant. Such an explanation was indeed found. 
In Theorem ~\ref{thmviro} we show how to express the encomplexed
writhe in terms of the shade number and the wrapping (self-linking)
number.   


\subsection*{Acknowledgments} The author would like to thank Oleg
Viro for useful suggestions, and John Etnyre and Jun Li for
stimulating discussions.


\section{Statements of main results}
In this section, some notation is introduced and the main results are
stated. 


\subsection{Linking of complex varieties with real projective
space}\label{IRempty} 
If $W$ is a projective $k$-dimensional variety in $(2k+1)$-dimensional
complex projective  
space then its set of closed points $\C W\subset\C P^{2k+1}$, with
the subspace topology and the induced complex orientation, is an oriented
$2k$-cycle in $\C P^{2k+1}$. Assume that 
$\R W=\C W\cap\R P^{2k+1}=\emptyset$ and let $[\C W]$ denote the
homology class of $\C W$ in $H_{2k}(\C P^{2k+1}-\R P^{2k+1})$.  
\begin{dfn}\label{dfnCshadeno}
The number
$$
\sh(W)=\frac12\left([\Gamma]\bullet[\C W]\right)\in\frac12\Z,
$$
where $[\Gamma]$ is the shade class and $\bullet$ denotes the
intersection product, is called  the {\em shade number} of $W$.
\end{dfn}

The range of the shade number can be expressed in terms of degrees:

\begin{thm}\label{thmrange}
For projective $k$-dimensional varieties in complex projective
$(2k+1)$-space of degree $d$ without 
real points, the range of the shade number consists of all
half-integers between $-\frac12d^2$ and $\frac12d^2$ which are
congruent to $\frac12d$ modulo $1$.
\end{thm}
\noindent
Theorem ~\ref{thmrange} is proved in \S ~\ref{pfthmrange}. 

\begin{prp}\label{prpinterpol}
If $W_0$ and $W_1$  are projective $k$-dimensional
varieties in complex projective $(2k+1)$-space without real points and
if $A$ is a 
$(2k+1)$-chain in $\C P^{2k+1}$ such that $\partial A=\C W_1-\C W_0$
then  
\begin{equation*}
\sh(W_1)-\sh(W_0)=\R P^{2k+1}\bullet A,
\end{equation*} 
where $\bullet$ denotes intersection number.
\end{prp}
\noindent
Proposition ~\ref{prpinterpol} is proved in \S ~\ref{pfprpinterpol}.

\begin{rmk}\label{rmkinterpol}
Proposition ~\ref{prpinterpol} gives information about how
the shade number changes under deformations of a variety.
If there is an ambient isotopy
$\phi_t$ of $\C P^{2k+1}$
carrying $\C W_0$ to $\C W_1$ in $\C P^{2k+1}$ 
(e.g. $\phi_t$ could be a $1$-parameter family of complex projective
transformations, or induced by a rigid isotopy) then $A$ can be taken 
as $A=\bigcup_t\phi_t(\C W_0)$. In particular, the shade number
changes by $\pm 1$ when a deformation crosses  
$\R P^{2k+1}$ transversely. So, in a sense, it is a first order
Vassiliev invariant (this notion is borrowed from knot
theory, see ~\cite{V}) of $k$-dimensional varieties in 
$\C P^{2k+1}$ without real points. 
\end{rmk}


\subsection{Shade numbers of armed real varieties}\label{Ishadearm}
If $V$ is a real projective subvariety of real projective $n$-space, then
after base extension, it is a projective subvariety of complex projective
$n$-space. We shall use the notions $\C V\subset\C P^n$ for its
set of closed points with the topology induced from the complex
manifold $\C P^n$, and $\R V$ to denote $\C V\cap\R P^n$
(often without further mentioning of the base extension).

Following ~\cite{A}, a submanifold $M$ of a
manifold $Y$ which is equipped with a non-vanishing normal vector
field $n$ 
will be called an {\em armed} submanifold. We denote it $(M,n)$.
Analogously, we say that a real projective variety $V$ in real projective
$m$-space is {\em armed} if its real part $\R V$ is an 
armed submanifold of $\R P^{m}$. We denote it $(V,n)$, where $n$
refers to the vector field. Note that, for dimensional
reasons, the real part of any $k$-dimensional subvariety of 
$(2k+1)$-dimensional projective space admits a non-vanishing normal
vector field.  

For real projective $k$-dimensional varieties $V$ in projective
$(2k+1)$-space, the intersection  
$\R V=\C V\cap\R P^{2k+1}$ may be large 
and there are many ways of pushing $\C V$ off 
$\R P^{2k+1}$ in $\C P^{2k+1}$. For armed varieties $(V,n)$, the
normal vector field 
$n$ determines such a push-off:  
$in$ ($i=\sqrt{-1}$) is a normal vector 
field of both $\C V$ and of $\R P^{2k+1}$
along $\R V$ in $\C P^{2k+1}$. 
If $\C V$, with its complex orientation, is shifted slightly along a
normal vector field $\nu$ 
extending $in$, and with support in a small neighborhood of $\R V$
in $\C V$, then an oriented 
$2k$-cycle $\C V_n$ in $\C P^{2k+1}-\R P^{2k+1}$ is obtained. 
Let $[\C V_n]$ denote its homology class in 
$H_{2k}(\C P^{2k+1}-\R P^{2k+1})$.  
\begin{dfn}\label{dfnRshadeno}
The number
$$
\sh(V,n)=\frac12\left([\Gamma]\bullet[\C V_n]\right)\in\frac12\Z,
$$
where $[\Gamma]$ is the shade class and $\bullet$ denotes the
intersection product, is called
{\em the shade number} of $(V,n)$.  
\end{dfn}
It is easy to see that $\sh(V,n)$ is independent of the choice of
$\nu$ (the extension of $in$) as long as its support and the shifting
distance are sufficiently small. 


\subsection{Real algebraic spheres}\label{IRsphere}
Spheres frequently occurs as real algebraic
ambient spaces  (e.g.\ the link of a complex plane curve singularity
is a real subvariety of the sphere rather than of projective space)
and so we define shade numbers of armed real varieties also in this
setting: the $n$-dimensional sphere $S^n$ is the set of real points of a
real quadric $Q^n$ in projective $(n+1)$-space (see \S
~\ref{shadeofRS}, for explicit equations), $S^n\approx\R Q^n$. This
real variety $Q^n$ will be called simply {\em the real
$n$-sphere}. Its set of complex 
points is the complex manifold $\C Q^n\subset\C P^{n+1}$. 

Using a projection $\C Q^{2k+1}\to\C P^{2k+1}$ which is a double
cover branched over a purely imaginary quadric and which restricts to
the standard double cover\linebreak 
$\R Q^{2k+1}\to\R P^{2k+1}$, 
we construct the
complex shade of $\R Q^{2k+1}$ in $\C Q^{2k+1}$ see Definition
~\ref{dfnshadeS} and Remark ~\ref{rmkshadePS}. 
Definition ~\ref{dfnRshadeno} then applies also to 
armed $k$-dimensional varieties in $\R Q^{2k+1}$. 

Our notation for subvarieties of the sphere is analogous to the one in
the projective case: if $V$ is a subvariety of the real $n$-sphere then 
$\C V\subset\C Q^n$ denotes its set of closed points (after base
extension) with topology induced from the complex manifold 
$\C Q^n$ and $\R V=\C V\cap\R Q^n$.

One may study also other real algebraic ambient spaces. 
In this paper we will however restrict attention to the two basic
cases of projective spaces and spheres.

\subsection{Even-dimensional real varieties}\label{Ieven}
If $(M,n)$ is an armed $2j$-dimensional submanifold of an oriented manifold
$Y$ of dimension $4j+1$ then let $e(M,n)$ denote the Euler number of the  
$2j$-dimensional vector bundle $\xi(n)$ which is the normal bundle of 
$M$ in $Y$ divided by its $1$-dimensional
subbundle generated by $n$. (Note that the orientation of $Y$ together 
with $n$ induce an orientation on the total space of $\xi(n)$, hence
the Euler number is defined, see \S ~\ref{eulno}.)  

\begin{thm}\label{thmeven}
Let $(V,n)$ be an armed $2j$-dimensional projective variety
without real 
singularities in real projective $(4j+1)$-space or in the real
$(4j+1)$-sphere. Then
\begin{equation*}
\sh(V,n)=(-1)^j\tfrac12e(\R V,n)
\end{equation*}
\end{thm}  
\noindent
Theorem ~\ref{thmeven} is proved in \S ~\ref{pfthmeven}.


\subsection{Odd-dimensional real varieties}\label{Iodd}  
For even-dimensional armed varieties, the shade number equals a topological
invariant. This is not the case in odd-dimensions. As we shall see,
for orientable odd-dimensional armed varieties, a linear
combination of the shade number and a topological invariant (the
wrapping number) is independent of the normal vector field
(which always exists, see \S ~\ref{Ishadearm}), and is invariant
under a certain class of deformations (weak rigid isotopies). 

The  {\em wrapping number} of an armed orientable $(2j+1)$-dimensional
submanifold $(M,n)$ of $S^{4j+3}$ ($\R P^{4j+3}$) is   
the orientation independent part of its self-linking
number. It will be denoted $\wra(M,n)$. In $S^{4j+3}$,
$\wra(M,n)\in\Z$, in $\R P^{4j+3}$,
$\wra(M,n)\in\frac12\Z$, see \S ~\ref{wrano}. With this
term introduced we define the linear combination mentioned above:

\begin{dfn}\label{dfnEW}
Let $(V,n)$ be an armed projective $(2j+1)$-dimensional
variety without real singularities in real projective
$(4j+3)$-space or in the real $(4j+3)$-sphere, with $\R V$
orientable. Define 
\begin{equation*}
\EW(V)=
\wra(\R V,n)+(-1)^{j}\sh(V,n)\in\tfrac12\Z.
\end{equation*}
\end{dfn}

A {\em weak rigid isotopy} of a real subvariety $V$ in a nonsingular
projective real algebraic variety $Y$ is a continuous $1$-parameter
family of real subvarieties $V_t$, $0\le t\le 1$, such that $V_0=V$,
and such that   
the induced $1$-parameter family $\R V_t\subset\R Y$ is given by
$\phi_t(\R V)$, where $\phi_t\colon\R Y\to\R Y$, $0\le t\le 1$, is a
$1$-parameter family of diffeomorphisms starting at the identity. 

For example, rigid isotopies are weak rigid isotopies but not conversely:
the $1$-parameter family $(\C V_t,\R V_t)$ induced by a rigid isotopy
satisfies $(\C V_t,\R V_t)=\psi_t(\C V,\R V)$, where $\psi_t$ is a
continuous $1$-parameter family of diffeomorphisms of the pair 
$(\C Y,\R Y)$. Also, continuous paths of real projective
transformations starting at the identity induce weak rigid isotopies in
projective space.   

\begin{thm}\label{thmodd}
Let $V$ be a projective $(2j+1)$-dimensional
variety without real singularities in real projective $(4j+3)$-space
or in the real $(4j+3)$-sphere, with $\R V$ orientable.
Then $\EW(V)$ is an integer, independent of
the choice of normal vector field, independent of the 
choice of orientation, and invariant under weak rigid isotopy. 
\end{thm}
\noindent
Theorem ~\ref{thmodd} is proved in \S ~\ref{pfthmodd}. For explicit
computations of $\EW$ 
of real algebraic representatives of the unknot and the
trefoil knot in $\R Q^3$, see 
\S ~\ref{Oknot} and \S ~\ref{32knot}, respectively. 

The invariant $\EW$ may change under deformations which are not
weak rigid isotopies. For certain deformations it changes in a controlled
manner. We begin by describing such deformations: 

Let $\epsilon>0$ and let $V_t$, $t\in(-\epsilon,\epsilon)$ be a
continuous $1$-parameter family of real
projective $(2j+1)$-dimensional varieties in  
real projective $(4j+3)$-space or in the real $(4j+3)$-sphere 
such that $V_t$ is without real singularities and $\R V_t$ is
connected and orientable, for $t\ne 0$. Assume that $V_0$ has
exactly one real 
double point where either two branches of $\R V_0$ intersect cleanly
or where two complex conjugate branches of $\C V_0$ intersect
cleanly. In the former case assume that the intersection point of the
traces of the deformations of the two branches  
of $\R V_t$ in $\R P^{4j+3}\times (-\epsilon,\epsilon)$ 
($\R Q^{4j+3}\times (-\epsilon,\epsilon)$)
is transverse, 
in the later assume that the intersection of the trace of (either) one
of the branches of $\C V_t$ meets 
$\R P^{4j+3}\times(-\epsilon,\epsilon)$ 
($\R Q^{4j+3}\times(-\epsilon,\epsilon)$)
transversely in 
$\C P^{4j+3}\times(-\epsilon,\epsilon)$
($\C Q^{4j+3}\times(-\epsilon,\epsilon)$).   
 
\begin{thm}\label{thmjump}
For $\epsilon>t>0$ and $V_t$ as described above
$$
\EW(V_t)-\EW(V_{-t})=\pm 2.
$$ 
\end{thm}  
\noindent
Theorem ~\ref{thmjump} is proved in \S ~\ref{pfthmjump}. 

\begin{rmk}
Theorem ~\ref{thmjump} implies in particular that $\EW$ is a first
order Vassiliev invariant of real algebraic embeddings $\phi$ 
of a real projective $(2j+1)$-dimensional variety $V$, with $\R V$
connected and 
orientable, into real projective $(4j+3)$-space, where $\phi$ varies
in a given linear system 
$\Lambda$ of dimension $N>4j+3$ on $V$ 
($\Lambda$ presents $V$ as a real algebraic variety without real
singularities in real projective $N$-space). 

This can be seen as follows: in $\Lambda$ there is a
{\em discriminant hypersurface} consisting of all maps $\phi$ 
with real singularities. The
top-dimensional strata of this discriminant consists 
of maps with one real-real or one complex-complex-conjugate double
point and $1$-parameter families intersecting the discriminant
transversely satisfies the assumptions on the deformations in  Theorem
~\ref{thmjump}. Therefore the first jump
of $\EW$ is $\pm 2$, and the second jump is zero  
(see ~\cite{Ar} for the notion of jump). That is, $\EW$ is a
non-trivial first order Vassiliev invariant. 
\end{rmk}


\subsection{Shade number and encomplexed writhe}
Shade numbers provide an
intrinsic 3-dimensional formula, homological in nature, 
for Viro's invariant of real algebraic links mentioned in 
\S ~\ref{applral}:
 
\begin{thm}\label{thmviro}
Let $V$ be a real algebraic link in real projective $3$-space. Then $\EW(V)$  
equals the encomplexed writhe of $V$ as defined in ~\cite{V}.  
\end{thm} 
\noindent
Theorem ~\ref{thmviro} is proved in \S ~\ref{pfthm2}. 
\begin{rmk}
There are different possible definitions of the encomplexed writhe of
a real algebraic link $V$ such that $\R V$ has more than one connected
component, see ~\cite{V}, \S 1.4. We choose one such definition,
see Remark ~\ref{choice}.   
\end{rmk}
In ~\cite{V}, Section 3.3, it is mentioned that the diagrammatic
approach to the encomplexed writhe can be generalized to give rigid
isotopy invariants of non-singular orientable real projective 
$(2j+1)$-dimensional varieties in real projective
$(4j+3)$-space. Carrying out the 
proposed generalization, invariants which agree with $\EW$
(Definition ~\ref{dfnEW}) are obtained. That is, the counterpart of Theorem
~\ref{thmviro} holds in the high-dimensional situation and two
formulas, one diagrammatic and one intrinsic homological, for the same
invariant are obtained.


\section{Complex shades and their characteristics}
In this section, orientation conventions are specified, complex shades
of projective spaces and of spheres are formally defined, and an
axiomatic characterization of the shade class is given. 


\subsection{Orientation conventions}\label{orconv}
Let $M$ be an oriented $k$-dimensional submanifold of an oriented
manifold $Y$ of dimension $n$. Let  
$NM$ denote the normal bundle of  $M\subset Y$, and let $x\in M$. A
basis $(n_1,\dots,n_{n-k})$  
of $N_xM$ is {\em positively oriented} if for any positively oriented
basis $(t_1,\dots,t_k)$ of the tangent space $T_xM$ of $M$ at $x$,
$(n_1,\dots,n_{n-k},t_1,\dots,t_k)$ is a positively 
oriented basis of $T_xY$.

If $M$ is an oriented manifold with boundary $\partial M$, we induce
an orientation on $\partial M$ by requiring that the {\em outward}
normal vector field of $\partial M$ in $M$ induces the positive
orientation of the normal bundle of $\partial M$ in $M$.


\subsection{Definition of the shade of projective space}\label{shadeofRP}
The complex shade of\linebreak 
$\R P^{2k+1}$ is constructed as follows.
Fix a point $p\in\R P^{2k+1}$.  
Let $L(p)\approx\R P^{2k}$ denote the set of all real lines through $p$.

If $l\in L(p)$ then $(\C l,\R l)\approx(\C P^1,\R P^1)$. Let
$X=\bigcup_{l\in L(p)}(\C l-\R l)$. Then $X$ is a disk-bundle over
the $2k$-sphere with fibers open 2-disks which are naturally
identified with connected components of $\C l-\R l$. Let $\bar X$ denote the
corresponding disk-bundle with fibers closed 2-disks. The fibers of
$\bar X$ are then naturally identified with the closures of components
of $\C l-\R l$ and there is a canonical map 
\begin{equation}\label{gammaP}
\gamma_p\colon\bar X\to\C P^{2k+1}.
\end{equation}
The restriction of $\gamma_p$ to $\pa\bar X-\gamma^{-1}_p(p)$
is a double cover of $\R P^{2k+1}-\{p\}$. The standard orientation 
of $\R P^{2k+1}$ therefore induces an orientation on $\partial\bar X$
which in turn induces an orientation on $\bar X$.  

\begin{dfn}\label{dfnshadeP}
The {\em complex shade} $\Gamma_p$ of $\R P^{2k+1}$, constructed using 
the point $p$, is the relative $(2k+2)$-cycle in 
$(\C P^{2k+1}, \R P^{2k+1})$ defined by   
$\Gamma_p=\gamma_p(\bar X)$, where $\gamma_p$ and $\bar X$ are as in
~\eqref{gammaP}. 
\end{dfn}

\begin{rmk}\label{rmkshadeP}
It is straightforward to check that $\Gamma_p$ in Definition
~\ref{dfnshadeP} has the following properties.
\begin{itemize}  
\item The homology class  of $\Gamma_p$ in 
$H_{2k+2}(\C P^{2k+1},\R P^{2k+1})$ is independent of $p$.  

\item The boundary $\partial \Gamma_p$ equals $2\R P^{2k+1}$.

\item If ${}^\ast\colon\C P^{2k+1}\to\C P^{2k+1}$ denotes complex
conjugation then $\Gamma^\ast_p=\Gamma_p$.
\end{itemize}
\end{rmk}


\subsection{Definition of the shade of the sphere}\label{shadeofRS}
Before we construct complex shades of spheres we introduce some
notation which will be used throughout the rest of the paper:

Let $Q^n$ denote the quadric in real projective $(n+1)$-space defined
by the homogeneous equation $-x_0^2+x_1^2+\dots+x_{n+1}=0$ in projective
coordinates $[x_0,\dots,x_{n+1}]$. Then $\R Q^n\approx S^n$.  

Let $\Pi\colon\C Q^n\to\C P^n$ denote the projection from
the point $[1,0,\dots,0]$ to the hyperplane in $\C P^{n+1}$
given by the equation $x_0=0$. Note that $\Pi$ is a double cover branched
over the purely imaginary quadric given by the equations
$x_1^2+\dots+ x^2_{n+1}=0, x_0=0$. If 
$\pi\colon\R Q^n\to\R P^n$ denotes the restriction of $\Pi$ then
$\pi$ is the standard double cover identifying antipodal points on 
$\R Q^n$. Thus, inverse images of lines in $\R P^n$ are great circles in
$\R Q^n$. 

We construct the shade of $\R Q^{2k+1}$. Fix a point $p\in\R P^{2k+1}$.
Let $\pi^{-1}(p)=\{\tilde p_0,\tilde p_1\}\subset\R Q^{2k+1}$. Then
$\tilde p_0$ and $\tilde p_1$ are antipodal points.  
Let $G(p)\approx\R P^{2k}$ denote the set of all great circles through
$\tilde p_0$ and $\tilde p_1$.

If $g\in G(p)$ then $g$ is the 
intersection of $\Delta$ and $Q^{2k+1}$, where
$\Delta$ is a real projective $2$-plane in $\R P^{2k+2}$ through
$\tilde p_0$, $\tilde p_1$, and one other point in $\R Q^{2k+1}$. Thus,
$g$ is a quadric in $\Delta$ and $(\C g,\R g)\approx(\C P^1,\R P^1)$.

Let $Y=\bigcup_{g\in G(p)}(\C g-\R g)$. Then $Y$ is a disk-bundle over
the $2k$-sphere with fibers open 2-disks which are naturally identified
with connected components of $\C g-\R g$. Let $\bar Y$ denote the
corresponding disk-bundle with fibers closed 2-disks. The fibers of
$\bar Y$ are then naturally identified with the closures of components
of $\C g-\R g$ and there is a canonical map 
\begin{equation}\label{gammaS}
\gamma_p\colon \bar Y\to\C Q^{2k+1}.
\end{equation}
The restriction of $\gamma_p$ to 
$\pa\bar Y-\gamma_p^{-1}(\{\tilde p_0,\tilde p_1\})$ 
is a (trivial) double cover of of 
$\R Q^{2k+1}-\{\tilde p_0,\tilde p_1\}$. 
The standard orientation of  
$\R Q^{2k+1}$ therefore induces an orientation on $\partial\bar Y$
which in turn induces an orientation on $\bar Y$.  

\begin{dfn}\label{dfnshadeS}
The {\em complex shade} $\Gamma_p$ of $\R Q^{2k+1}$, constructed using 
$p\in\R P^{2k+1}$, is the relative
$(2k+2)$-cycle in $(\C Q^{2k+1}, \R Q^{2k+1})$ defined by  
$\Gamma_p=\gamma_p(\bar Y)$, where $\gamma_p$ and $\bar Y$ are as in
~\eqref{gammaS}.  
\end{dfn}

\begin{rmk}\label{rmkshadeS}
Remark ~\ref{rmkshadeP} carries over word by word from the projective to the
spherical case if ``$\R P^{2k+1}$''  and ``$\C P^{2k+1}$'' are replaced
by ``$\R Q^{2k+1}$'' and ``$\C Q^{2k+1}$'', respectively.
\end{rmk}

\begin{rmk}\label{rmkshadePS}
Let $p\in\R P^{2k+1}$ and let $\Gamma_p^{\po}$ and
$\Gamma_p^{\qu}$ denote the
shades of $\R P^{2k+1}$ and $\R Q^{2k+1}$,
respectively, constructed using $p$. Then 
$\Gamma_p^{\qu}=\Pi^{-1}(\Gamma_p^{\po})$. 
Moreover, the map $\Pi$ restricted to 
any $\C g\subset\Gamma_p^{\qu}$ is the standard double cover of
a $\C l\subset\Gamma_p^{\po}$, branched at two distinct
complex conjugate points.
\end{rmk}


\subsection{Homology characterization of the shade}
Let $L^\po$ denote the complex $k$-dimensional subspace in 
complex projective $(2k+1)$-space given by the equations  
\begin{equation}
iz_0-z_1=iz_2-z_3=\dots=iz_{2k}-z_{2k+1}=0,
\end{equation}
in projective coordinates $[z_0,\dots,z_{2k+1}]$. 
Then $\R L^\qu=\emptyset$.  

Let $L^\qu$ denote the complex $k$-dimensional subspace in 
complex projective $(2k+2)$-space given by the equations  
\begin{equation}
iz_1-z_2=iz_3-z_4=\dots=iz_{2k+1}-z_{2k+2}=0,\quad z_0=0.
\end{equation}
in projective coordinates $[z_0,\dots,z_{2k+2}]$.  
Then $\C L^\qu\subset\C Q^{2k+1}$, see \S
~\ref{shadeofRS}, and $\R L^\qu=\emptyset$. 

Let $[L^\po]$ and $[L^\qu]$ 
denote the homology classes in 
$H_{2k}(\C P^{2k+1}-\R P^{2k+1})$ and 
$H_{2k}(\C Q^{2k+1}-\R Q^{2k+1})$ of
$\C L^\po$ and $\C L^\qu$, with their complex orientations,
respectively. Let $[F^\po]\in H_{2k}(\C P^{2k+1}-\R P^{2k+1})$
and $[F^\qu]\in H_{2k}(\C Q^{2k+1}-\R Q^{2k+1})$
denote the homology classes
corresponding to the fiber classes in the normal bundles of 
$\R P^{2k+1}$ and $\R Q^{2k+1}$ in $\C P^{2k+1}$ and $\C Q^{2k+1}$,
respectively.   

To prove the next lemma we need the homology of $\C Q^{2k+1}$.
A straightforward calculation 
(the adjunction formula gives the total
Chern class of  
$\C Q^{2k+1}$ and, in particular, its Euler characteristic $2k+2$,
apply Lefschetz hyperplane theorem, Poincar{\'e} duality and the
universal coefficient theorem) 
shows that
\begin{equation}\label{homCS}
H_r(\C Q^{2k+1})=\begin{cases}
                 \Z  & \text{if $r$ is even,}\\
                  0  & \text{if $r$ is odd.}
                 \end{cases}
\end{equation}

\begin{lma}\label{lmaHC-R}
Let $n=2k+1$. The homology groups $H_{n-1}(\C P^n-\R P^n)$ and
$H_{n-1}(\C Q^n-\R Q^n)$ are isomorphic to $\Z\oplus\Z$, are
generated by $[L^\qu]$ and $[F^\qu]$ respectively $[L^\qu]$ and
$[F^\qu]$, and the homomorphisms  
\begin{align*}
\alpha^\po\colon H_{n+1}(\C P^n,\R P^n)\to
\Hom(H_{n-1}(\C P^n-\R P^n);\Z);&\quad
\alpha^\po(\Sigma)(\Delta)=\Sigma\bullet\Delta,\\ 
\alpha^\qu\colon H_{n+1}(\C Q^n,\R Q^n)\to
\Hom(H_{n-1}(\C Q^n-\R Q^n);\Z);&\quad
\alpha^\qu(\Sigma)(\Delta)=\Sigma\bullet\Delta,\\ 
\end{align*}
where $\bullet$ denotes the intersection product, are isomorphisms.
\end{lma}

\begin{pf}
The two cases will be treated simultaneously. Therefore, we let
$(\C Y^n,\R Y^n)$ denote 
either $(\C P^n,\R P^n)$ or $(\C Q^n,\R Q^n)$, and we drop the
superscripts on $F$, $L$, and $\alpha$.

Poincar\'e duality implies that the following diagram with exact
rows commutes and that all vertical arrows (which are cap-products
with the orientation class of $\C Y^n$) are isomorphisms.  
$$
\begin{CD}
\dots H^{2n-r-1}(\R Y^n) @>>{\delta}> H^{2n-r}(\C Y^n,\R Y^n) @>>>
H^{2n-r}(\C Y^n)\dots \\ 
{} @VVV @VVV @VVV {}\\ 
\dots H_{r+1}(\C Y^n,\C Y^n-\R Y^n) @>>{\partial}> 
H_{r}(\C Y^n-\R Y^n) @>>>
H_{r}(\C Y^n) \dots 
\end{CD}.
$$ 
It follows that $H_{n-1}(\C Y^n-\R Y^n)\approx\Z\oplus\Z$, generated as
claimed, and that $H_{n-2}(\C Y^n-\R Y^n)=0$, see Equation
~\eqref{homCS}.  
The universal coefficient theorem then implies that 
$$
\Hom(H_{n-1}(\C Y^n-\R Y^n);\Z)\approx 
H^{n-1}(\C Y^n-\R Y^n).
$$
Let $M^{2n}$ be the compact manifold which is the complement of an
open tubular neighborhood of $\R Y^n$ in $\C Y^n$. By homotopy,
$$
H^{\ast}(M^{2n})\approx H^{\ast}(\C Y^n-\R Y^n).
$$
Poincar{\'e} duality gives 
$H^{r}(M^{2n})\approx H_{2n-r}(M^{2n},\partial M^{2n})$. The inclusion 
$M^{2n}\to\C Y^n$ induces a map 
$(M^{2n},\partial M^{2n})\to (\C Y^n,T)$, where $T$ is a closed tubular
neighborhood of $\R Y^n$ in $\C Y^n$. By excision,
$$
H_{\ast}(M^{2n},\partial M^{2n})\approx H_{\ast}(\C Y^n,T).
$$
It now follows from homotopy and the 5-lemma that  
$$
H_{\ast}(\C Y^n, T)\approx H_{\ast}(\C Y^n,\R Y^n).
$$
Hence,
$$
H_{n+1}(\C Y^n,\R Y^n)\stackrel{\alpha}{\approx}
\Hom(H_{n-1}(\C Y^n-\R Y^n);\Z). 
$$ 
\end{pf}

Remark ~\ref{rmkshadeP} (Remark ~\ref{rmkshadeS}) implies that the
homology class of the shade $\Gamma_p$ is independent of $p$. We
denote this homology class $[\Gamma]$.  
If $[\R P^{2k+1}]\in H_{2k+1}(\R P^{2k+1})$ ($[\R Q^{2k+1}]\in
H_{2k+1}(\R Q^{2k+1})$) denotes the orientation class and ${}^\ast$
denotes complex conjugation then this class can be characterized as
follows.

\begin{prp}\label{prpHax}
The {\em shade class} $[\Gamma]$ is the unique class in 
$H_{2k+2}(\C P^{2k+1},\linebreak
\R P^{2k+1})$ ($H_{2k+2}(\C Q^{2k+1},\R Q^{2k+1})$) 
which satisfies the following conditions 
\begin{itemize}
\item[{\rm (a) }] $[\Gamma]=[\Gamma]^\ast$,
\item[{\rm (b) }] $\partial[\Gamma]=2[\R P^{2k+1}]\quad 
\left(\partial[\Gamma]=2[\R Q^{2k+1}]\right)$,
\item[{\rm (c) }] $\left|[\Gamma]\bullet[L^\po]\right|=1
\quad\left(\left|[\Gamma]\bullet[L^\qu]\right|=1\right)$,
\end{itemize}
where $\bullet$ denotes the intersection product.
\end{prp}

\begin{rmk}\label{rmkHax}
If $k$ is even then (c) is a
consequence of (a) and (b). Note that Proposition ~\ref{prpHax}
could be taken as definition of the shade class.   
\end{rmk}

\begin{pf}
We use notation as in the proof of Lemma ~\ref{lmaHC-R}. Also, let
$J^\po$ be a real $(k+1)$-dimensional linear subspace of real
projective $(2k+1)$-space, and  $J^\qu$ be the intersection of $Q^{2k+1}$ 
and the $(k+2)$-dimensional linear subspace, of 
real projective $(2k+2)$-space, given by the equations
$x_2=x_4=\dots=x_{2k}=0$ in projective coordinates
$[x_0,\dots,x_{2k+2}]$. We write $[J]$ for the homology class 
represented by $\C J^\po$ ($\C J^\qu$) with its complex orientation in
$H_{n+1}(\C Y^{n},\R Y^n)$.

By Remarks ~\ref{rmkshadeP} and ~\ref{rmkshadeS},
$\partial[\Gamma]=2[\R Y^{2k+1}]$. Hence $[\Gamma]$ satisfies (a).  

Let $n=2k+1$ and consider the exact sequence
$$
0\to H_{n+1}(\C Y^{n}) \to
H_{n+1}(\C Y^{n},\R Y^{n})\stackrel{\partial}{\to}H_n(\R Y^n)\to 0.
$$
Since $\partial[\Gamma]=2[\R Y^{n}]$, any homology class $\xi$ such
that $\partial\xi=2[\R Y^n]$ can be written as 
\begin{equation}\label{kernel}
[\Gamma]+a[J], 
\end{equation}
where $a\in\Z$. (In the spherical case, note that
$\C J^\qu\bullet\C L^\qu=1$, thus $[J]$ generates the image of
$H_{n+1}(\C Q^n)$.) 
 
By Proposition ~\ref{lmaHC-R}, the group $H_{n-1}(\C Y^n-\R Y^n)$ dual to 
$H_{n+1}(\C Y^{n},\R Y^{n})$ is generated by $[F]$ and $[L]$. It is
straightforward to check that 
$[\Gamma]\bullet[F]=[\Gamma]^\ast\bullet[F]$ and that 
$[\Gamma]\bullet[L]=[\Gamma]^\ast\bullet[L]$. It follows from this
that $[\Gamma]=[\Gamma]^\ast$. Hence $[\Gamma]$ satisfies condition
(b). Moreover, $[J]^\ast=(-1)^{k+1}[J]$. Thus,
Equation ~\eqref{kernel} implies that the homology class $[\Gamma]$ is
uniquely determined by (a) and (b) if $k$ is even.
  
Finally, 
$\left|[\Gamma]\bullet[L]\right|=1$. Hence $[\Gamma]$ satisfies
condition (c). Equation ~\eqref{kernel} together with the fact
that $[J]\bullet[L]=1$ 
imply that $[\Gamma]$ is uniquely determined by (a)-(c),
if $k$ is odd.
\end{pf}

\begin{rmk}\label{Hconj}
The proof of Proposition ~\ref{prpHax} gives a partial picture of the
relation between conjugation and $2k$-dimensional homology. The
following two equations complete it:
\begin{align*}
[F]^\ast &=-[F],\\
[L]^\ast &=\begin{cases}
              -[L] & \text{if $k$ is odd},\\ 
              {[L]}-([\Gamma]\bullet[L])[F] & \text{if $k$ is even},
             \end{cases}
\end{align*}
The first equation is obvious. To see that the second equation holds
we consider first projective space and
use a $(2k+1)$-cycle in $\C P^{2k+1}$ interpolating between 
$\C L^\po$ and 
$\C (L^\po)^\ast$, where $\C (L^\po)^\ast$ denotes the complex linear
subspace conjugate to $L$, with its complex orientation.

More precisely, consider the sequence of maps 
$\phi_j\colon\C P^{k}\times [-1,1]\to\C P^{2k+1}$ defined as
follows. For $0\le j\le k$, let $p_j$ denote the coordinate vector
in coordinates 
$[z_0,\dots,z_{2k+1}]$ on $\C P^{2k+1}$ which has $z_{2j}=1$,
$z_{2j+1}=i$, and 
$z_m=0$ for $m\ne 2j,2j+1$, let $p_j^\ast$ be the complex conjugate
coordinate vector, and let $q_j$ denote the coordinate vector which has
$z_{2j+1}=i$ and $z_m=0$ for $m\ne 2j+1$. Let $[u]=[u_0,\dots,u_k]$ be
homogeneous coordinates on $\C P^k$ and define 
\begin{equation*}
\phi_j([u],t)=
u_0 p_0^\ast +\dots +u_{j-1} p_{j-1}^\ast + 
u_j (p_j-(1+t)q_j) + u_{j+1}p_{j+1} +\dots + u_{k} p_k,
\end{equation*} 
where the left-hand side is to be interpreted as a coordinate vector.

Note that $\phi_j([u],1)=\phi_{j+1}([u],-1)$. The cycle 
$\phi_j(\C P^{k}\times [-1,1])$ intersects $\R P^{2k+1}$ 
transversely in exactly one point and  the sign of the
intersection points of $\phi_j$ and that of $\phi_{j+1}$ are opposite. 
Finally, $\phi_0(\C P^k,-1)$ is $\C L^\po$ and $\phi_k(\C P^{k},1)$ is
$\C (L^\po)^\ast$. Since $[L^\po]^\ast$ equals $(-1)^k$ times the homology
class of $\C (L^\po)^\ast$ it follows that 
$[L^\po]^\ast=-[L^\po]$ if $k$ is odd and, since
$\left|[L^\po]\bullet[\Gamma]\right|=1$ and $[\Gamma]\bullet[F^\po]=2$, 
$[L^\po]^\ast=[L^\po]-([\Gamma]\bullet[L^\po])[F^\po]$ if $k$ 
is even, as claimed.

The argument in $\C Q^{2k+1}$ is similar: write the coordinates on 
$\C P^{2k+2}$ as $[w,z]$ where $z=[z_0,\dots,z_{2k+1}]$ 
is as above, and use the maps
$$
\psi_j([u],t)=\left[\sqrt{1-t^2},\phi_j([u],t)\right],
$$
where $\phi_j([u],t)$ is as above, to interpolate between $\C L^\qu$
and $\C (L^\qu)^\ast$.
\end{rmk}


\section{Varieties without real points}
In this section, Proposition ~\ref{prpinterpol}  and
Theorem ~\ref{thmrange} are proved.


\subsection{Proof of Proposition ~\ref{prpinterpol}}\label{pfprpinterpol}
We use the notation $[F^\po]$ as in Lemma ~\ref{lmaHC-R}.  
If the intersection number of the $(2k+1)$-chain $A$ with 
$\R P^{2k+1}$ is $m$ then $[\C W_1]-[\C W_0]=m[F^\po]$ in 
$H_{2k}(\C P^{2k+1}-\R P^{2k+1})$. Since $[\Gamma]\bullet[F^\po]=2$
the statement follows.\qed


\subsection{Proof of Theorem ~\ref{thmrange}}\label{pfthmrange}
Let $L^\po$ be as in Lemma ~\ref{lmaHC-R}. Then $\sh(L^\po)=\pm\frac12$.
Let $W$ be a $k$-dimensional variety of degree $d$ with 
$\R W=\emptyset$. By 
definition of degree, $\C W$ is homologous 
to $d\cdot\C L^\po$ in $\C P^{2k+1}$.  Proposition
~\ref{prpinterpol} then implies that $\sh(W)$ is congruent to
$\frac12 d$ modulo $1$. 

Let $p\in\R P^{2k+1}$ be a point and $H$ a hyperplane in 
real projective $(2k+1)$-space such that $p\notin\R H$. 
Let $\pr\colon(\C P^{2k+1}-\{p\})\to\C H$ denote linear projection
and let $W^\ast$ denote the complex conjugate variety of $W$. 

Possibly after a small algebraic perturbation, we may assume that
$\pr(\C W)$ meets 
$\pr(\C W^\ast)$ transversely. 
(Proposition ~\ref{prpinterpol} implies that small perturbations
does not affect $\sh(W)$.) 

If $x$ is an intersection point of $\Gamma_p$ and $\C W$ then $x$
lies on $\C l$ for some real line $l$ through $p$. The complex
conjugate $x^\ast$ of $x$ satisfies $x^\ast\in\C l$ and 
$x^\ast\in \C W^\ast$. Now, $\pr(\C l)$ is a point and hence
$\pr(x)=\pr(x^\ast)$. Thus, to each intersection point of $\Gamma_p$
and $\C W$ there corresponds an intersection point of $\pr(\C W)$
and $\pr(\C W^\ast)$ in $\C H$. The $k$-dimensional varieties $\pr(W)$ and
$\pr(W^\ast)$ both have degree $d$. They meet transversely, and hence
intersect in $d^2$ points. It follows that $|\sh(W)|\le\frac12 d^2$.

To see that the shade number takes all possible values between
$-\frac12 d^2$ and 
$\frac12 d^2$ we construct explicit varieties: 

Let $K$ be a large positive real number. Define 
\begin{equation*}
P(u,v)=K(u-v)(u-2v)\dots(u-dv).  
\end{equation*}
For $t\in\R$, define 
\begin{equation*}
Q_t(u,v)=K\left(u-(t+\tfrac{1}{d^2+1})v\right)
\left(u-(t+\tfrac{2}{d^2+1})v\right)\dots
\left(u-(t+\tfrac{d}{d^2+1})v\right).
\end{equation*}
Then the equations $P(u,1)=1$ and $Q_t(u,1)=1$ has $d$ real solutions
$\theta_1<\dots<\theta_d$ and $\phi_1(t),\dots,\phi_d(t)$,
respectively. The distances $|\theta_j-j|$ and 
$\left|\phi_j(t)-\left(t+\frac{j}{d^2+1}\right)\right|$,
$j=1,\dots,d$, can be made arbitrarily small by choosing $K$
sufficiently large. 

For $t\in\R$, let $W_t$ be the variety in complex projective
$(2k+1)$-space defined by the equations
\begin{align*}
& P(z_1,z_0) + iQ_t(z_2,z_0) -(1+i)z_0^d = 0,\\
& z_1 + z_2 - iz_3 =0,\\
&z_{2j}-iz_{2j+1}=0,\quad\text{ for }j=2,\dots k.
\end{align*}
in projective coordinates 
$[z_0,\dots,z_{2k+1}]$ (with complex conjugation given by conjugation
on all coordinates).

Consider the shade $\Gamma_p$ constructed using the point
$p=[0,0,0,1,0,\dots,0]$ in $\R P^{2k+1}$. A
straightforward check shows that $\Gamma_p\cap \C W_t$ consists of the
points $[1,\theta_j,\phi_k(t),i(\theta_j+\phi_k(t)),0,\dots,0]$, 
$1\le j,k\le d$. For large negative $t$, all intersection points
contribute with the same sign to $\Gamma_p\bullet \C W_t$ and
$\sh(W_t)=\pm\frac12 d^2$. As $t$
increases, the roots $\phi_1(t),\dots,\phi_d(t)$ approaches $-\theta_d$.
At some $t'$, $\phi_d(t')=-\theta_d$ and $W_{t'}$ has a real point 
$[1,\theta_d,\phi_d(t'),0,0,\dots,0]$. As $t$ passes $t'$ the sign of the
intersection point
$[1,\theta_d,\phi_d(t),i(\theta_d+\phi_d(t)),0,\dots,0]$ changes. This 
changes the intersection number $\Gamma_p\bullet \C W_t$ by $\mp 2$ and
hence $\sh(W_t)$ by  
$\mp 1$. As $t$ increases further the roots $\phi_{d-1}(t)$,\dots,
$\phi_1(t)$ passes $-\theta_d$ and at each passage $\sh(W_t)$
changes by $\mp 1$. As $t$ increases even further the same things
happen at $-\theta_{d-1}$, \dots, $-\theta_1$ until $\phi_1(t)$ has
passed $-\theta_1$ and we have $\sh(W_t)=\mp\frac12d^2$.\qed 


\section{Topological invariants, shade numbers, twists, and generic shades}
In this section, two topological invariants of
armed submanifolds are described. For later use, the behavior of
these invariants 
and of shade numbers under local deformations of normal vector fields are 
studied. 

The existence of shades that have good properties with respect
to a given real variety is established and used to demonstrate a
certain symmetry property of shade numbers in even dimensions.


\subsection{Euler numbers}\label{eulno}
For the readers convenience we recall the definition of Euler
number: let $M$ be a $k$-dimensional manifold and let $\xi$ be
$k$-dimensional vector bundle over $M$. Assume that the total space
$E(\xi)$ of $\xi$ is orientable. 

Let $\mathcal M$ and $\mathcal F$ denote the integer local coefficient
systems over $M$ associated to the 
orientation bundle of $M$ and the fiber orientation bundle of $\xi$,
respectively. Then the orientation of $E(\xi)$ gives 
\begin{equation}\label{eqnlocsys}
{\mathcal M}\otimes{\mathcal F}\approx\Z,
\end{equation} 
where $\Z$ denotes the
(trivialized) local coefficient system associated to the the
orientation bundle of $E(\xi)$, restricted to the zero-section $M$.  
An isomorphism ${\mathcal M}\approx{\mathcal F}$, is specified by
requiring that, at each point, ~\eqref{eqnlocsys} 
is given by ordinary multiplication of integers. 

This specified isomorphism in turn gives a well-defined pairing
$$
\la\, ,\ra\colon H^k(M;{\mathcal F})\otimes 
H_k(M;{\mathcal M})\to\Z. 
$$ 
The {\em Euler number} $e(\xi)$ of $\xi$ is defined as 
$\la e,[M]\ra$, where $e$ is the Euler class of the bundle (the
obstruction to finding a non-vanishing section), and
$[M]$ is the orientation class of $M$.

One can compute the Euler number of $\xi$ by choosing any 
section $s$ of $\xi$ transverse to the zero-section and sum up the
local intersection numbers at zeros of $s$. The local intersection
number at a zero $p$ of $s$ is the intersection number in $E(\xi)$ of
some neighborhood $U\subset M$ of $p$, with any chosen orientation,
and $s(U)$, with the orientation induced from the chosen orientation 
on $U$.

It is straightforward to check that if $k$ is odd then $e(\xi)=0$.   


\subsection{Wrapping numbers}\label{wrano}
Let $(M,n)$ be an armed $k$-dimensional orientable\linebreak 
smooth
submanifold of $\R P^{2k+1}$ or $S^{2k+1}$.
Let $K(1)$, \dots, $K(m)$ denote the connected components of $M$
and let $K(j)_n$ denote a copy of 
$K(j)$ shifted slightly along $n$. Fix some orientation of
$M$ and note that it induces an orientation on each $K(j)_n$.  

\begin{dfn}\label{dfnfrno}
The number
\begin{equation*}
\wra(M,n)=\sum_{j=1}^m\lk(K(j),K(j)_n).
\end{equation*}
is called the {\em wrapping number} of $(M,n)$. 
\end{dfn}

\begin{rmk}\label{rmkwraor}
Note that $\wra(M,n)$ is independent of the choice of
orientation on $M$. If $M$ is
endowed with an orientation then other self-linking numbers can be  
associated to it. They are expressible in terms of
$\wra(M,n)$ and the linking numbers $\lk(K(j),K(l))$, $j\ne l$.
\end{rmk}

\begin{rmk}\label{rmkwraeul}
An argument similar to the proof of Theorem ~\ref{thmeven}
(\S ~\ref{pfthmeven}) shows that if $\dim(M)$ is even then
$\wra(M,n)=-\tfrac12e(M,n)$. Note that 
the definition of the wrapping number applies only if $M$ is
orientable whereas the Euler number is defined also for non-orientable 
manifolds. 
\end{rmk}


\subsection{Twists}
Let $M$ be a $k$-dimensional submanifold of an oriented manifold
$Y$ of dimension $2k+1$. Let 
$NM$ denote the normal bundle of $M\subset Y$, and let 
$N^0M$ denote the bundle of  non-zero vectors in $N M$.

If $n$ and $m$ are sections in $N^0M$ then a standard obstruction
theory argument shows that $n$ can be made homotopic to $m$ by local
changes which we call local twists: 

Pick a Riemannian metric in the normal bundle $NM$ and let 
$UN M$ denote the corresponding 
bundle of unit normal vectors. Let $v$ be a section of $UN M$ and let
$p\in M$. Choose a local coordinate neighborhood $X\subset M$
around $p$ and a  
local trivialization $\beta\colon X\times S^k\to UNM|X$ such that 
$\beta^{-1}(v(x))=(x,a)\in X\times S^{k}$ for $x$ varying in $X$ and
$a$ fixed in $S^{k}$.   
Fix an orientation on $X$. This orientation together with the
orientation on $Y$ induce an orientation on the fibers of
$NM|X$, which in turn induces an orientation on the fibers $S^k$ of
$UNM|X$.  
Let $D^{k}\subset X$ be a disk centered at $p$ and let $r\colon D^k\to
S^k$ be a map of degree $\pm 1$ such that $r(\partial D^k)=a$. Define
$v'$ as 
$$
v'(x)=\begin{cases}
        v(x), & x\in M-D^k\\
        \beta(x,r(x)), & x\in D^k
      \end{cases}.
$$
We say that
$v'$  {\em is obtained from $v$ by adding a local twist at $p$} and
that the degree of $r$ is {\em the sign of the local twist}.

We shall now study how Euler numbers, wrapping numbers and shade
numbers are affected by local twists. 
  
\begin{lma}\label{lmaeultw}
Let $M$ be a $2j$-dimensional submanifold of an oriented manifold
$Y$ of dimension $4j+1$. Let 
$m$ and $n$ be non-vanishing normal vector fields of $M$. 
\begin{itemize}
\item[{\rm (a)}] If $m$
and $n$ are homotopic as sections in $N^0M$ then
$e(M,m)=e(M,n)$. 
\item[{\rm (b)}]
If $m$ is obtained from $n$ by adding a 
local twist of sign $\sigma=\pm1$ then $e(M,m)=e(M,n)+2\sigma$.
\item[{\rm (c)}] $e(M,n)=-e(M,-n)$.
\end{itemize}
\end{lma}
\begin{pf}
Assertion (a) is immediate. To prove (b), let $p$ be the point where
the twist is added and fix local coordinates
$x=(x',x'',x_{4j+1})\in\R^{4j+1}$, $x'=(x_1,\dots,x_{2j})$,
$x''=(x_{2j+1},\dots,x_{4j})$, on an open set $X\subset Y$ such 
that $M\cap X=\{x\colon x''=0,x_{4j+1}=0\}$, and $p$ corresponds
to $x=0$. 

Write $\pa_i=\tfrac{\pa}{\pa x_i}$. We may assume that $n$ is given by 
$x'\mapsto \frac{1}{\sqrt{2j}}\sum_{k=2j+1}^{4j}\pa_k$ in 
coordinates $x'$ on $M\cap X$. 
Identify the bundle $\xi(n)$ with the orthogonal
complement of $n(x')$ in the normal bundle and choose a section $s$ of
$\xi(n)$ which is given by $x'\mapsto\pa_{4j+1}$ in $M\cap X$.

Let $S^{2j}=\{x\colon x'=0, (x'')^2+x_{4j+1}^2=1\}$ and  let $D^{2j}$
be a small disk where the twist, $x'\mapsto r(x')\in S^{2j}$ is added
to $n$ to obtain $m$. The bundle $\xi(m)$ is then identical to
$\xi(n)$ outside $D^{2j}$ and over $D^{2j}$, $\xi(m)$ can be
identified with the bundle with fiber over points $x'$ equal to
$T_{r(x')}S^{2j}$. Let $s'$ be 
the section of $\xi(m)$ which agrees with $s$ outside $D^{2j}$, and
inside $D^{2j}$ is given by orthogonal projection of $s(x')$ into
$T_{r(x')}S^{2j}$. That is,
$$
s'(x')=\pa_{4j+1}-\left\la\pa_{4j+1},r(x')\right\ra r(x'),
$$
where $\la\, ,\ra$ denotes the Euclidean metric on $X$. 
Then $s'(x')=0$ if and only if $r(x')=\pm\pa_{4j+1}$. 
At such an $x'$, 
$$
ds'(x')=\mp dr(x').
$$
Since changing the sign of each vector in an even dimensional frame
does not change its orientation, it follows that the difference of the
algebraic number of zeroes of $s'$ and $s$ equals twice the 
degree of $r$. This implies (b).

Assertion (c) follows from the fact that the bundles 
$\xi(n)$ and $\xi(-n)$ are canonically isomorphic, but the orientations of
the total space $E(\xi(n))=E(\xi(-n))$ induced by $n$ and $-n$,
respectively, are opposite. 
\end{pf}

\begin{lma}\label{lmawratw}
Let $M$ be an orientable $(2j+1)$-dimensional submanifold of
$S^{4j+3}$ or $\R P^{4j+3}$. Let 
$m$ and $n$ be non-vanishing normal vector fields of $M$. 
\begin{itemize}
\item[{\rm (a)}] If $m$
and $n$ are homotopic as sections in $N^0M$ then
$\wra(M,m)=\wra(M,n)$. 
\item[{\rm (b)}]
If $m$ is obtained from $n$ by adding a 
local twist of sign $\sigma=\pm1$ then
$\wra(M,m)=\wra(M,n)+\sigma$. 
\end{itemize}
\end{lma}
\begin{pf}
Assertion (a) is immediate. To prove (b), let $p$ be the point where
the twist is added and fix local coordinates
$x=(x',x'',x_{4j+3})\in\R^{4j+3}$, $x'=(x_1,\dots,x_{2j+1})$,
$x''=(x_{2j+2},\dots,x_{4j+2})$, on an open set $X$ in the ambient space 
such that $M\cap X=\{x\colon x''=0,x_{4j+3}=0\}$, and $p$ corresponds
to $x=0$. 

Let the orientation of the ambient space along $X$
be given by the frame\linebreak 
$(\pa',\pa'',\pa_{4j+3})$, where $\pa'=(\pa_1,\dots,\pa_{2j+1})$ and
$\pa''=(\pa_{2j+2},\dots,\pa_{4j+2})$. We may assume that
$n$ is given by
$x'\mapsto\tfrac{1}{\sqrt{2j+1}}\sum_{k=2j+2}^{4j+2}\pa_k$ in
coordinates $x'$ on $M\cap X$. Choose the orientation on
$M$ given by the frame $\pa'$.  

Let $K$ denote the component of $M$ which intersects $X$. 
To find 
$\lk(K^{2j+1}, K^{2j+1}_n)$ we need a $(2j+2)$-chain $A$
with boundary $2K$. We may choose $A$ so that 
$A\cap X=a_1\cup a_2$,  where
\begin{align*}
a_1 &=\{x\colon x''=0,x_{4j+3}\ge 0\},\\
a_2 &=\{x\colon x''=0,x_{4j+3}\le 0\}.
\end{align*}
Then, by our orientation conventions (see \S ~\ref{orconv}),
$(\partial',\partial_{4j+3})$ and $(\partial',-\partial_{4j+3}$
are positively oriented frames along $a_1$ and $a_2$, respectively.

Let $S^{2j+1}=\{x\colon x'=0, (x'')^2+x_{4j+3}^2=1\}$ and let
$D^{2j+1}$ be a small disk where
the twist $x'\mapsto r(x')\in S^{2j+1}$ is added to $n$ to obtain
$m$. If $U$ is a small 
neighborhood of $D^{2j+1}$ in $X$ then outside of $U$, $K_n=K_m$,
and inside $U$,
$$
K_n=\{x\colon x''=\delta'',x_{2k+1}=0\},
$$ 
where $\delta''=(\delta,\dots,\delta)$ for some small $\delta>0$,
and
$$
K_m=\{x\colon (x'',x_{2k+1})=\delta r(x')\}.
$$
Thus $U\cap K_n\cap A=\emptyset$, and $U\cap K_m\cap A$ consists of
points $(x',\delta r(x'))$ such that
$r(x')=\pm\pa_{4j+3}$. Assume that $\pm\partial_{4j+3}$ are
regular values of $r$.

Consider a point $x'$ such that $r(x')=\partial_{4j+3}$. The
corresponding intersection point 
$(x',\delta r(x'))$ lies in $a_1$ and the local intersection
number of $A$ and $K_m$ is given by the 
sign of the orientation of the frame 
\begin{equation}\label{sign1}
(\pa',\pa_{4j+3},\pa'+dr\partial'),
\end{equation} 
which is just the local degree of $r$ at $x'$. At a point 
$(x',\delta r(x'))$ such that $r(x')=-\partial_{4j+3}$ the sign of
$\partial_{4j+3}$ in Formula ~\eqref{sign1} must be changed 
since the intersection point lies in $a_2$ and not $a_1$ and again the 
local intersection number equals the local degree of $r$.

It follows that the intersection numbers $A\bullet K_m$ and
$A\bullet K_n$ differs by twice the degree of $r$. Since there
is a factor $\frac12$ in the definition of $\lk$ (see Section ~\ref{intr}), 
(b) follows.   
\end{pf}

\begin{lma}\label{lmashadetw}
Let $(V,n)$ be a $k$-dimensional projective variety without real
singularities in real projective $(2k+1)$-space or in the real
$(2k+1)$-sphere. Write $k=2j+1$ if $k$ is odd and
$k=2j$ if $k$ is even. If $m$ is obtained by adding a local twist of
sign $\sigma=\pm 1$ to $n$ then
$$
\sh(V,m)=\sh(V,n)+(-1)^{k+j}\sigma.
$$
\end{lma}
\begin{pf}
Let $\R Y^{2k+1}$ and $\C Y^{2k+1}$ denote the real- and the complex
points of the ambient space, respectively. Let $p\in \R V$ be the
point where the twist is added and let $K$ be 
the connected component of $\R V$ which contains $p$.

Introduce
(real analytic) local coordinates $x=(x',x'',x_{2k+1})\in\R^{2k+1}$, 
where $x'=(x_1,\dots,x_k)$ and $x''=(x_{k+1},\dots,x_{2k})$, 
on an open set $X\subset\R Y^{2k+1}$ so that
$\R V\cap X=\{x\colon x''=0,x_{2k+1}=0\}$, and so that $p$
corresponds to $x=0$ (and hence $\R V\cap X=K\cap X$).  
Then there exist holomorphic local coordinates
$z=(z',z'',z_{2k+1})\in\C^{2k+1}$ 
on a neighborhood $Z\subset \C Y^{2k+1}$ of $p$ 
such that $Z\cap\R Y^{2k+1}=X$, such that $z_j=x_j+iy_j$,
$j=1,\dots,2k+1$, such that the complex conjugation on 
$\C Y^{2k+1}$ is given by conjugation of coordinates in $Z$, and such
that $\C V$ is locally given 
by $\C V\cap Z=\{z\colon z''=0,z_{2k+1}=0\}$.
We use the notions $\pa'$ and $\pa''$ as in the proof of
Lemma ~\ref{lmaeultw}. Let $(\pa',\pa'',\pa_{2k+1})$ be 
a positively oriented frame of $T\R Y^{2k+1}|X$.

The shade $\Gamma=\Gamma_q$ ($q\in\R P^{2k+1}$) of $\R Y^{2k+1}$ is a
relative $(2k+2)$-cycle in $(\C Y^{2k+1},\R Y^{2k+1})$. After a homotopy,
supported in a small 
neighborhood of $Z$, we may assume that $\Gamma\cap Z$ equals 
$\gamma_1\cup\gamma_2$, where 
\begin{align*}
\gamma_1 &= \{z=x+iy\colon  y'=y''=0,y_{2k+1}\ge 0 \},\\
\gamma_2 &= \{z=x+iy\colon y'=y''=0, y_{2k+1}\le 0 \}.
\end{align*}
Since
$\partial\Gamma=2\R Y^{2k+1}$, our orientation conventions
(\S ~\ref{orconv}) imply that
$$
(\partial',\partial'',\partial_{2k+1},i\partial_{2k+1})\text{ and }
(\partial',\partial'',\partial_{2k+1},-i\partial_{2k+1})
$$
are positively oriented frames of $\gamma_1$ and $\gamma_2$,
respectively. 

We may assume that the normal vector field 
$n$ is given by
$x'\mapsto \frac{1}{\sqrt{k}}\sum_{j=k+1}^{2k} \partial_{j}$ in
coordinates $x'$ on $K\cap X$. 
Let $S^k=\{x\colon x'=0, (x'')^2+x_{2k+1}^2=1\}$ and let
$D^k$ be a small disk where
the twist $x'\mapsto r(x')\in S^{2j+1}$ is added to $n$ to obtain
$m$. 

If $E$ is a small 
neighborhood of $D^{2j+1}$ in $Z$ then outside of $E$, $\C V_n=\C V_m$,
and inside $E$,
$$
\C V_n=\{z\colon x''=z_{2k+1}=0, y''=\delta''\}  
$$
and
$$  
\C V_m=
\{z\colon x''=x_{2k+1}=0, (y'',y_{2k+1})=i\phi(y')r(x')\},
$$
where $\phi$ is a positive function such that 
$\phi(y')=\delta$ for $(y')^2<\epsilon$ and $\phi(y')= 0$
for $(y')^2>2\epsilon$ for some very small $\epsilon>0$.
Thus,  $E\cap\C V_n\cap\Gamma=\emptyset$, 
and $E\cap\C V_m\cap \Gamma$
consists of points $(x'+0i,0+i\delta r(x'))$, where
$r(x')=\pm\pa_{2k+1}$. Hence, the difference between $\sh(V,m)$ and
$\sh(V,n)$ equals one half times the sum of local intersection
numbers of $\C V_m$ and $\Gamma$ in $E$. 
We calculate this difference: 

Let $dx'\wedge dy'=dx_1\wedge dy_1\wedge\dots\wedge dx_k\wedge dy_k$,
interpret  $dx''\wedge dy''$ and $dx\wedge dy$ similarly. Then the
complex orientation of 
$\C Y^{2k+1}$ is given by the $(4k+2)$-form 
$dx\wedge dy$
and the complex orientation of $\C V$ is given by $dx'\wedge dy'$.

Assume that $r(x')=\partial_{2k+1}$. Then the local intersection
number of $\Gamma$ and $\C V_m$ at the corresponding intersection
point is given by the sign of  
\begin{align}\label{sign2}
&(-1)^{j} dx\wedge dy
\left(\pa',\pa'',\pa_{2k+1},i\pa_{2k+1},\pa'+idr(x')\pa', i\pa'\right)\\
&=
(-1)^{j}
dx\wedge dy
\left(\pa',i\pa',\pa'',idr(x')\pa',\pa_{2k+1},i\pa_{2k+1}\right),\notag
\end{align}
where $(-1)^{j}$ (recall $k=2j$, if even, and $k=2j+1$, if odd)
arises since 
$dx'\wedge dy'=(-1)^{j}dx_1\wedge\dots\wedge dx_k\wedge
dy_1\wedge\dots\wedge dy_k$.    
The local degree of $r$ at $x'$ is given by the sign of
\begin{equation}\label{sign3}
dx_{k+1}\wedge\dots\wedge dx_{2k+1}\left(\pa_{2k+1},dr(x')\pa'\right)=
(-1)^k dx_{k+1}\wedge\dots\wedge dx_{2k}\left(dr(x')\pa'\right).
\end{equation}
It follows that the local intersection number of $\Gamma$ and 
$\C V_m$ at the intersection point corresponding to $x'$ equals
$(-1)^{j+k}$ times the local degree of $r$ at $x'$. 

At an intersection point with
$r(x')=-\partial_{2k+1}$ a 
similar calculation shows that local intersection number is again
$(-1)^{j+k}$ times the local degree of $r$ at $x'$. 
(Both the sign of $i\partial_{2k+1}$ in Formula ~\eqref{sign2} and
that of $\pa_{2k+1}$ in Formula \eqref{sign3} change.) 
\end{pf}


\subsection{Generic shades}\label{generic}
We shall treat projective spaces and spheres simultaneously. The
notation introduced in \S ~\ref{shadeofRS} will be used repeatedly.

Let $V$ be a real projective $k$-dimensional variety without real
singularities  
in real projective $(2k+1)$-space (in the real $(2k+1)$-sphere). We
define three subsets of $\R P^{2k+1}$, which we want to avoid. 

First, consider the variety $\Sigma$ of singular points of $V$. 
By assumption $\R\Sigma=\emptyset$. Moreover, the dimension of
$\Sigma$ is at most $k-1$. Let $A$ be the chordal variety (see
~\cite{H}) of $\Sigma$ ($\Pi(\Sigma)$). Then the complex dimension of 
$\C A$ is at most $2k-1$ and $\R A$ is a stratified set of real
dimension at most $2k-1$.  

Second, let $B$ be the tangential variety (see ~\cite{H}) of
$V$ ($\Pi(V)$). Then the
complex dimension of $\C B$ is at most $2k$ and 
$\R B$ is a stratified set of real dimension at most
$2k$.  

Third, let $D$ 
be the variety which is the
closure of all points on chords through two distinct smooth points 
$p,q\in\C V$ ($p,q\in\Pi(\C V)$) such that the projective tangent
spaces of $\C V$ ($\Pi(\C V)$) at $p$ and $q$ intersect.
Then the complex dimension of $\C D$ is at most $2k$ and 
$\R D$ is a stratified set of real dimension at most $2k$.  

A point $p\in\R P^{2k+1}$ will be called {\em generic} (with respect
to $V$) if $q\notin \R A\cup \R B\cup \R D$. Note that the set of generic 
points form a dense open set in $\R P^{2k+1}$.

If $p\in\R P^{2k+1}$ is a generic point 
then the shade $\Gamma_p$ 
of $\R P^{2k+1}$ ($\R Q^{2k+1}$), constructed
using $p$, does not intersect $\C\Sigma$ (since $p\notin \R A$), and at
any  
point $x\in\Gamma_p\cap(\C V-\R V)$ 
($x\in\Gamma_p\cap (\C V-\R V)$)
the intersection is transverse (since $p\notin \R D$).  

Assume that $\R V$ is orientable and oriented, and let $K$ be a
component of $\R V$. To describe linking properties of 
$K$ we need a chain with boundary $2K$, see Section
~\ref{intr}. Such chains can be constructed using generic points.

First consider $\R P^{2k+1}$: 
if $y,z\in\R P^{2k+1}$ are distinct then let $l(y,z)$ denote the line
containing them. Let $C=\bigcup_{x\in K}\R l(p,x)$.
If $p$ is a generic point then $C$ is an immersed submanifold, except
possibly at $p$ (since $p\notin \R B$). Clearly, we 
can orient the chain $C$ so that $\partial C=2K$.  

Second consider $\R Q^{2k+1}$:
if $y\in\R P^{2k+1}$, $z\in\R Q^{2k+1}$, and $z\notin\pi^{-1}(y)$,
then let $g(y,z)$ denote the great circle through $\pi^{-1}(y)$ and $z$.
Let $p\in\R P^{2k+1}$ and let 
$C=\bigcup_{x\in K}\R g(p,x)$. If $p$ is generic
then $C$ is an immersed submanifold, except possibly at $\tilde p_0$ and
$\tilde p_1$ (since $p\notin \R B$). Let $\cle^1$ be the unit circle in
$\C$. Let 
$\phi\colon K\times\cle^1\to C$ be a map such that $\phi(x,1)=x$,
$\phi(x,i)=\tilde p_0$, $\phi(x,-i)=\tilde p_1$, and such that
$\phi(x,-)\colon \cle^1\to \R g(p,x)$ is an embedding, for every
$x$. Consider $M=K\times\cle^1-(K\times\{i\})$. Then $M$ is a
manifold and $K^{k}\times\{1\}$ subdivides $M$ into two components
$M'$ and $M''$. Orient these so that they both induce the positive
orientation on their boundary $K\times\{1\}$. Endow the chain $C$
with the orientation induced from $M'$ and $M''$. Then $C$ is a chain
with $\partial C=2K$. 

\begin{rmk}\label{rmknearfar}
If $(V,n)$ is an armed real projective $k$-dimensional variety without
real 
singularities in $(2k+1)$-dimensional real projective space or in the
real $(2k+1)$-sphere and if $p\in\R P^{2k+1}$ is a 
generic point, then there are two types of intersection points in
$\Gamma_p\cap \C V_n$ and in $C\cap \R V_n$: points near $\R V$
and points far from $\R V$. More precisely, since $p$ is generic
there is a small neighborhood $E$ of $\R V$ in $\C P^{2k+1}$ 
($\C Q^{2k+1}$) such that $\Gamma_p\cap\C V\cap E=\R V$. We write 
$\Gamma_p(E)=\Gamma_p\cap E$ and $\Gamma'_p(E)=\Gamma_p-\Gamma_p(E)$
and say that a point $x\in\Gamma_p\cap\C V_n$ is of type 
\begin{itemize}
\item[(\cf)] if $x\in\Gamma'(E)\cap\C V $, and of type
\item[(\cn)] if $x\in\Gamma(E)\cap\C V$.  
\end{itemize}

We use similar notation for the chain $C$: let $U=E\cap\R P^{2k+1}$
($U=E\cap\R Q^{2k+1}$) and 
write $C(U)=C\cap U$ and $C'(U)=C-C(U)$ and say that a point 
$y\in C\cap K_n$ is  of type
\begin{itemize}
\item[(\rf)] if $x\in C'(U)\cap K_n$, and of type
\item[(\rn)] if $x\in C(U)\cap K_n$. 
\end{itemize}
\end{rmk}

\subsection{Symmetry in even dimensions}
\begin{lma}\label{lmaevensym}
Let $(V,n)$ be an armed projective $2j$-dimensional variety without real
singularities in real projective $(4j+1)$-space or in the real
$(4j+1)$-sphere. Then 
$$
\sh(V,n)=-\sh(V,-n).
$$  
\end{lma}
\begin{pf}
Pick a generic point $p$ in $\R P^{4j+1}$ and let $\Gamma=\Gamma_p$ be
the shade constructed using $p$. Let $\R Y^{4j+1}$ and $\C Y^{4j+1}$
denote the real- and complex parts of the ambient space,
respectively. Let $E$ be a neighborhood of 
$\R V$ in $\C Y^{4j+1}$ with properties as in 
Remark ~\ref{rmknearfar}. 

Assume that the vector field $\nu$ (which extends $in$) used to
shift $\C V$ has support inside $\C V\cap E$. 
As in Remark ~\ref{rmknearfar}, we distinguish two types (\cf) and (\cn)
of intersection points in $\C V\cap\Gamma$: 

If $x$ is an intersection point of type (\cf) then, since both $\Gamma$ 
and $\C V$ are invariant under conjugation, also the complex
conjugate point $x^\ast$ of $x$ is an intersection point of type
(\cf). We compare their signs: let $o(\Gamma,x)$ and $o(\C V,x)$ denote
positively oriented frames of $T_x\Gamma$ and $T_x\C V$,
respectively. The intersection number at $x$ is then given by the sign 
$\sigma$ of the frame 
$\left(o(\Gamma,x),o(\C V,x)\right)$. Since the complex
dimension of $\C Y^{4j+1}$ is odd,  
complex conjugation ${}^\ast$ reverses orientation and hence the sign
of the frame  $\left(o(\Gamma_p,x)^\ast,o(\C V,x)^\ast\right)$ at 
$x^\ast$ is $-\sigma$. The orientations of $o(\Gamma,x^\ast)$ and
$o(\Gamma,x)^\ast$ agree (recall, $[\Gamma]$ is invariant under
conjugation). Since $\C V$ is complex even-dimensional the
orientations $o(\C V,x)^\ast$ and $o(\C V,x^\ast)$ agree as
well. It follows that points of type (\cf) does not contribute at all to 
$\sh(V,n)$.

Consider an intersection point $q$ of type (\cn). Let $Z$ be a small
neighborhood of $q$ and let $X=\R Y^{4j+1}\cap Z$. We use coordinates
on $X$ and $Z$ as in the proof of Lemma ~\ref{lmashadetw}, more
precisely:  

Let $\R V\cap X=\{x\colon x''=0,x_{4j+3}=0\}$ in coordinates
$x=(x',x'',x_{4j+1})$ on $X$, and let
be $(z',z'',z_{4j+1})$, where $z'=x'+iy'$, $z'=x'+iy'$, and
$z_{4j+1}=x_{4j+1}+iy_{4j+1}$ be holomorphic coordinates on $Z$. Also, 
$\Gamma\cap Z=\gamma_1\cup\gamma_2$, where 
\begin{align*}
\gamma_1 &= \{z=x+iy\colon  y'=y''=0,y_{4j+1}\ge 0 \},\\
\gamma_2 &= \{z=x+iy\colon y'=y''=0, y_{4j+1}\le 0 \}.
\end{align*}

We view the normal vector field $n$ as a map 
$x'\mapsto n(x')\in\{x\colon x'=0, (x'')^2+x_{4j+1}^2=1\}$. It is then 
clear that the intersection point $q$ corresponds to a point $x'$ such 
that $n(x')=\pm\pa_{4j+1}$. Assume that $n(x')=\pa_{4j+1}$ and that
the intersection is transverse. The contribution from $q$ to
$\Gamma\bullet\C V_n$ is then given by the sign of
$$
(-1)^j dx\wedge dy
(\pa',\pa'',\pa_{4j+1},i\pa_{4j+1},idn\pa',i\pa').
$$  

Now consider the vector field $-n$. It is easy to see that to each
intersection point in $\C V_n\cap\Gamma$ there is a
corresponding intersection point of the same type, (\cf) or (\cn), in 
$\C V_{-n}\cap\Gamma$. As we have seen above, only the points
of type (\cn) contribute. 

Consider the intersection point $q'$ corresponding to the point $q$ (of type
(\cn) considered above). The contribution from $q'$ to 
$\Gamma\bullet\C V_{-n}$ is given by the sign of 
$$
(-1)^j dx\wedge dy
(\pa',\pa'',\pa_{4j+1},-i\pa_{4j+1},-idn\pa',i\pa'),
$$  
where the sign in front of $i\pa_{4j+1}$ appears since
$q'\in\gamma_2$ (in contrast to $q$, which lies in $\gamma_1$). 
Since multiplying all 
vectors in an even-dimensional frame by $-1$ does not change its
orientation, the contributions from $q$ and $q'$ have opposite signs. 
The lemma follows.  
\end{pf}


\section{Real varieties}

In this section, Theorems ~\ref{thmeven}, ~\ref{thmodd}, and
~\ref{thmjump} are proved. The diagrammatic definition of the
encomplexed writhe is presented, following  ~\cite{V} and, with this
definition at hand, Theorem ~\ref{thmviro} is proved.  

\subsection{Proof of Theorem ~\ref{thmeven}}\label{pfthmeven}
Let $\mathcal N$ denote the set of homotopy classes of sections in the
bundle of non-zero normal vectors of $\R V$. 

Let $\Phi,\Psi\colon {\mathcal N}\to\tfrac12\Z$ be defined as
\begin{align*}
\Phi({\mathbf n}) &=\tfrac12 e(\R V,n),\\
\Psi({\mathbf n}) &=\sh(V,n),
\end{align*}
where $n$ is a representative of the homotopy class $\mathbf n$. Then 
it follows from  Lemma ~\ref{lmaeultw} (a) that $\Phi$ is well-defined 
and from Lemma ~\ref{lmawratw} (a) that $\Psi$ is
well-defined. Moreover, Lemma ~\ref{lmaeultw} (b) and Lemma
~\ref{lmawratw} (b) imply that $\Phi-(-1)^j\Psi$ is a constant
function equal to $c$ say. But then Lemma ~\ref{lmaeultw} (c) and
Lemma ~\ref{lmaevensym} implies that 
\begin{align*}
c &=\tfrac12e(\R V,n)-(-1)^j\sh(V,n)=
\tfrac12e(\R V,-n)-(-1)^j\sh(V,-n)\\
 &= -(\tfrac12e(\R V,n)-(-1)^j\sh(V,n))=-c.
\end{align*}
Hence, $c=0$ and the theorem follows.\qed

\subsection{Proof of Theorem ~\ref{thmodd}}\label{pfthmodd}
It follows from Lemmas ~\ref{lmawratw} and
~\ref{lmashadetw} that $\EW$ is independent of the choice of normal
vector field of $\R V$. Also, $\EW$ is easily seen to be
independent of orientation on $\R V$, and to be invariant under 
weak rigid isotopy. 

It remains to show that $\EW$ is integer-valued: 
fix a generic point $p\in\R P^{4j+3}$ and let $\Gamma=\Gamma_p$ denote
the shade constructed using $p$. Let
$E$ be a small neighborhood of $\R V$ in $\C Y^{4j+3}$ as in
Remark ~\ref{rmknearfar}. As there, two 
types (\cf) and (\cn) of points in $\Gamma\cap\C V_n$ will be
distinguished. As in the proof of Lemma ~\ref{lmaevensym}, the points
of type (\cf) come in pairs, the same 
argument as there shows that, in the odd-dimensional case under
consideration, the two points in a pair of type (\cf) contribute with
the same sign to $\Gamma\bullet\C V_n$.

Let $K$ be a component of $\R V$. To calculate 
$\lk(K,K_n)$ we consider $C\bullet K_n$, where
$C$ is constructed using the generic point $p$, see \S
~\ref{generic}. Let $U=E\cap\R Y^{4j+3}$. As in Remark
~\ref{rmknearfar} there are two types, (\rf) and (\rn), of intersection
points in $C\cap K_n$ (we will use the notions $C(U)$,
$C'(U)$, $\Gamma(E)$, and $\Gamma'(E)$ as there).

We consider first points of type (\rf). Let $x'$ be such a point.  
\begin{itemize}
\item If ambient space is projective space then $x'$ lies very close to a
chord $l$ of $K$, which passes through $p$.
\item If the ambient space is the sphere then $x'$ lies very close to a
great circle $g$ which passes through two points on $K$ and
through $\tilde p_0$ and $\tilde p_1$ (see \S ~\ref{shadeofRS} for
notation). 
\end{itemize}

Let $x$ and $y$ be the points in $K$ on $\R l$ ($\R g$) and let
$x$ be the one close to $x'$. Then there is another intersection point 
$y'$ of type (\rf) close to $y$. Thus, points of type (\rf) come in pairs. 

We now consider points of types (\cn) and (\rn). Using the local models
of Lemmas ~\ref{lmawratw} and ~\ref{lmashadetw}, we observe that to each
intersection point $z$ of type (\rn) in 
$C(U)\cap K_n$ there 
corresponds exactly one point $r$ of type (\cn) in 
$\Gamma(E)\cap\C V_n$. It follows that $\EW$ is integer-valued. 
\qed

\begin{rmk}\label{rmksigns}
For future reference we establish the relation between the local
intersection numbers in the pairs considered in the above proof:

Let $x'$ and $y'$ be a pair of intersection points 
of type (\rf), as in the proof above. Then the local intersection
numbers of $C'(U)$ and $K_n$ at $x'$ and $y'$ equal the local
intersection numbers of $C'(U)$ and $K$ at $x$ and $y$,
respectively. These agree and therefore {\em the intersection points
contribute with the same sign to the intersection number 
$C\bullet K_n$}  

To see that this is the case, pick a sub-arc
$a$ of $\R l$ 
($\R g$) between $x$ and $y$, which does 
not contain $p$ ($\tilde p_0$, see \S ~\ref{generic}). Let $X$
and $Y$ be positively oriented frames of $T_xK$ and  $T_yK$,
respectively.    

\begin{itemize}
\item If the ambient space is projective space, $a$
lies in an affine 
part of $\R P^{4j+1}$ which does not contain $p$. If $u$ is a
non-vanishing tangent vector field of $a$ then the local
intersection number of $C$ and $K$ at $x$ is given by the
orientation of the frame $(\pm u,Y,X)$ and the local intersection number
at $y$ is given by the frame $(\mp u,X,Y)$. Since $K$ is
odd-dimensional, these frames have the same orientation. 
\item If the ambient space is the sphere,
we use stereographic projection\linebreak 
$s\colon \R Q^{2k+1}\to\R^{2k+1}$ from
$\tilde p_0$ such that $s(\tilde p_1)=0$. Two separate cases must be
considered.  

First, assume that $\tilde p_1\notin a$. If $u$ is a non-vanishing
tangent vector field of $s(a)$ then the orientation of the frame
$(\pm u,ds(Y),ds(X))$ gives the local intersection number of $C$ and
$K$ at $x$ and the orientation of the frame $(\mp
u,ds(X),ds(Y))$ gives the local intersection number at $y$.

Second, assume that $\tilde p_1\in a$. If $u$ is a non-vanishing
tangent vector field of $s(a)$ then the orientation of the frame
$(\pm u,-ds(Y),ds(X))$, where $-ds(Y)$ indicates the frame which is
obtained if all vectors in $ds(Y)$ are multiplied by $-1$,
gives the local intersection number of $C$ and
$K$ at $x$ and the orientation of the frame 
$(\mp u,-ds(X), ds(Y))$ gives the local intersection number at $y$. 

In any case, since $K$ is odd-dimensional, the local
intersection numbers are of the same sign.  
\end{itemize}

Consider intersection points $z$ of type (\cn) and $r$ of type (\rn), as
in the end of the above proof. Using the local models of Lemmas
~\ref{lmawratw} and ~\ref{lmashadetw} (cf. also the Proof of Theorem
~\ref{thmviro}) one calculates that the local intersection number of 
$C(U)$ and $\R V_n$  at $r$ equals $\lambda$ if and only if
the local intersection number of $\Gamma$ and $\C V_n$ at $z$
equals $(-1)^{j+2j+1}\lambda=(-1)^{j+1}\lambda$. Thus, the total
contribution of a pair of a point of type (\cn) and a corresponding
point of type (\rn) to $\EW$ equals $0$.  
\end{rmk}

\subsection{Proof of Theorem ~\ref{thmjump}}\label{pfthmjump}
We must check how the invariant $\EW(V_t)$ changes at the
double point instance. We separate the cases:
\begin{itemize}
\item[({\sc rr})] $V_0$ has a {\em real-real double point}, where 
two real branches pass through each other, or
\item[({\sc cc})] $V_0$ has a {\em complex-complex-conjugate double
point},  where two 
non-singular complex conjugate branches pass through each other at a
point in $\R P^{4j+3}$. 
\end{itemize} 

We use the coordinate notation as in Lemma ~\ref{lmawratw}. 
There exist (analytic) local coordinates
$x+i\,y=z=(z_1,\dots,z_{4j+3})=(z',z'',z_{4j+3})\in\C^{4j+3}$ 
in a neighborhood $Z\subset\C P^{4j+3}$
of the singular point $s\in \R V_0$ with the following
properties in cases (a) and (b) above.
  
\begin{itemize}
\item[({\sc rr})] The two local branches of $\R V_t$ near $s$
are given by the 
equations $x''=0,x_{4j+3}=t$ and $x'=0,x_{4j+3}=-t$, respectively.
\pagebreak

\item[({\sc cc})] The two local complex conjugate branches $B_t$ and
$B^\ast_t$ of  
$\C V_t$ near $s$ are given by the equations
$z'+iz''=0,z_{4j+3}=it$ and $z'-iz''=0,z_{4j+3}=-it$, respectively. 
\end{itemize}

Moreover, after changing the relative cycle $\Gamma$ by homotopy we
may assume that $\Gamma\cap Z=\gamma_1\cup \gamma_2$,
where 
\begin{align*}
\gamma_1&=\{z=x+iy\colon y'=y''=0,y_{4j+3}\ge 0\},\\
\gamma_2&=\{z=x+iy\colon y'=y''=0,y_{4j+3}\le 0\}
\end{align*}

In case (a), choose the normal vector field $n$ so that
$n=\pa_{2j+2}$ along 
the first local branch and $n=\pa_{1}$ along the second one. 
At $t=0$, the two local branches experiences a crossing change. 
Since $\R V_t$ is connected it follows that 
$\wra(\R V_t,n)=\lk(\R V_t, (\R V_t)_n)$ changes by
$\pm 2$.  
Clearly, $\sh(V_t,n)$ is 
constant in $t$ and thus, $\EW(V_t)$ changes by $\pm 2$ at $t=0$ in
case ({\sc rr}).  

In case ({\sc cc}) we may assume that the normal vector field $\nu$ of 
$\C V_t$, which is the extension of the vector field $i n$
along $\R V_t$, vanishes along $B_t$ and $B^\ast_t$. Note that, for
$t\ne 0$, there are two complex conjugate points  
$b_t$ and $b^\ast_t$ in 
$\Gamma\cap (\C V_t)_n\cap Z=\Gamma\cap\C V_t\cap Z$, 
$b_t\in\gamma_1$ and $b^\ast_t\in\gamma_2$. 

The orientations of $\gamma_1$ and $\gamma_2$ are ``opposite'' (see
the proof of Lemma ~\ref{lmashadetw}) and since the
dimension of $\C V_t$ is odd, complex conjugation changes
orientation. Therefore the local intersection numbers at $b_t$ and
$b^\ast_t$ agree. 

As $t\to 0+$, $b_t$ and $b^\ast_t$ come closer together. At $t=0$, the two
points collide at $s$, and as $t$ becomes negative $b_t$ has moved to
$\gamma_2$ and $b^\ast_t$ to $\gamma_1$. It follows that
$2\sh(V_t,n)=[\Gamma]\bullet[(\C V_t)_n]$ changes by
$\pm 4$ at $t=0$. Clearly,  
$\wra(\R V_t,n)$ is constant in $t$. 
Thus, $\EW$ changes by $\pm 2$ at $t=0$ also in case ({\sc cc}).
\qed


\subsection{The diagrammatic definition of encomplexed
writhe}\label{diagapp}  
We present the definition of the local writhe at a double point in a
projection of a real algebraic knot, see ~\cite{V}, Sections 2.1 and
2.2.  

Let $K$ be a variety in real projective $3$-space such that $\R K$ is
a knot  and let $c$ be a generic point in $\R P^3$. Fix some  hyperplane
$H$ in $\R P^3$, $c\notin \R H$. Let $p_c\colon(\R P^3-\{c\})\to \R H$ and 
$P_c\colon (\C P^3-\{c\})\to \C H$ denote the projections.
The image of $K$ under projection may have two kinds of double points:
\begin{itemize}
\item[(\rdp)] double points $x\in \R H$ which are transverse intersection of
two branches of $p_c(\R K)$, and  
\item[(\sdp)] double points $s\in \R H$ which are intersections of two
non-singular complex
conjugate branches of $P_c(\C K)$, both meeting $\R H$ transversely in
$\C H$ (such double points are called {\em solitary}).   
\end{itemize}

Let $x$ be a double point of type (\rdp) 
and let $l$ be the line which is the preimage
of $x$ under the projection. Denote by $a$ and $b$ the points of
$\R l\cap \R K$. The 
points $a$ and $b$ divide $\R l$ into two segments. Choose one of these
and denote it $S$. Choose an orientation of $\R K$. Let $v$ be $w$ be
nonzero tangent vectors of $\R K$ at $a$ and $b$ respectively, which are
directed along the selected orientation of $\R K$. Let $u$ be the tangent
vector of $\R l$ at $a$ directed inside $S$. Let $w'$ be a vector at $a$
which is tangent to the plane containing $l$ and $w$ and directed to
the same side of $S$ as $w$ (in the affine part of the plane
containing $S$ and $w$). The local writhe at $x$ is the sign of the
orientation of the 
frame $(v,u,w')$ in $T_a\R P^3$. It is independent of the choices of
$a$ and $S$, and of the choice of orientation of the knot.

Let $s\in \R H$ be a double point of type (\sdp)
and let $l$ be the preimage of $s$ under the projection. Then $\C l$ is the
preimage of $s$ under $P_c$ and $\C l\cap\C K$
consists of two imaginary complex conjugate points. Call them $a$ and
$b$. Since $a$ and $b$ are conjugate they belong to different components of 
$\C l-\R l$. Let $\D l(a)$ and $\D l(b)$ denote the closures of the
components of $\C l-\R l$ containing $a$ and $b$ respectively.
The complex orientations of $\D l(a)$ and $\D l(b)$ induce orientations
$u_a$ and $u_b$ respectively on $l$ (their common boundary). 
Let $\C K(a)$ and $\C K(b)$ denote neighborhoods in $\C K$ of $a$ and $b$,
respectively. The
images $P_c(\C K(a))$ and $P_c(\C K(b))$ intersects $\R H$
transversely in $\C H$ at $s$. Let $o_a$ be the local 
orientation of 
$\R H$ at $s$ which gives a positive local intersection number of
$\R H$ and $P_c(\C K(a))$ in $\C H$ (with its complex orientation). Define
$o_b$ similarly using $P_c(\C K(b))$ instead. The signs of
the orientations $(o_a,u_a)$ and $(o_b,u_b)$ of $\R P^3$ at $s$
agree. This sign is the local writhe at $s$.

It turns out that the sum of local writhes of all crossing points in a
generic projection is independent of the projection chosen, see
~\cite{V}, Sections 2.3 and 3.1. This sum is the {\em encomplexed
writhe} of the knot. 

\begin{rmk}\label{choice}
For $V$ such that $\R V$ is a many component link, the local
writhe at a double point of type (\sdp) is still well-defined. The local
writhe at a double point of type (\rdp) is well-defined only if the
preimages of the double point belong to the same component of $\R V$. 
To incorporate other double points in the definition of encomplexed
writhe 
one must fix a semi-orientation (an orientation up to sign) on
$\R V$, see \cite{V}, Section 1.4. However, the local writhes then depend
on the semi-orientation.   

Another possibility is to neglect these other double points and define
the encomplexed writhe of a real algebraic link as the sum of local
writhes at double points of type (\sdp) and at those double points of
type (\rdp) for which the 
preimages belong to the same component of $\R V$. {\em This is the
definition we use in the present paper.}  It has the advantage of
depending only on the algebraic link itself, not on chosen
orientations.  
\end{rmk}


\subsection{Proof of Theorem ~\ref{thmviro}}\label{pfthm2}
Fix an orientation of $\R V$ and endow $\R P^3$ with the metric coming
from the standard metric on $S^3\approx\R Q^3$. Choose a generic point
$c\in \R P^3$ and a hyperplane $H$, $c\notin \R H$. As in \S
~\ref{diagapp}, the notions $p_c$ and $P_c$ denote the real and
complex projections from $c$ to $\R H$ and $\C H$, respectively.

Let $K$ be a connected component of $\R V$. We shall define a
normal vector field $n$ along
$K$. Consider $C=\bigcup_{k\in K}l(c,k)$ as in \S ~\ref{generic}.  
Let $NK$ denote the bundle of vectors in $T\R P^3|K$ orthogonal to
the tangent vector of $K$. Let $\nu K$ denote the bundle of vectors in
$T\R P^3|K$ tangent to $C$ and orthogonal to the tangent vector of
$K$. Since $c$ is generic, $\nu K$ is a 1-dimensional sub-bundle of
the 2-dimensional bundle $NK$. Let $\xi K$ denote the orthogonal
complement of $\nu K$ in $NK$. It is easy to see that the bundles 
$\xi K$ and $\nu K$ are trivial if and only if $K$ is contractible in 
$\R P^3$.  

We make separate choices of $n$ for contractible and non-contractible
components $K$ of $\R V$. 
\begin{itemize}
\item If $K$ is contractible then let $n$ be a non-zero section in
$\xi K$.
\item If $K$ is not contractible then fix a small sub-arc $J$ of $K$
between $k_1$ and $k_2$, say, and an orientation preserving
diffeomorphism $[0,\pi]\to J$, $0\mapsto k_1$. The restriction of $n$
to $K-J$ is defined to be the (unique up to sign) non-zero section of 
$\xi K|(K-J)$. Choose a trivialization $NK|J\approx[0,\pi]\times\R^2$
such that if $(e_1,e_2)$ is the standard basis of $\R^2$ then $\R e_1$
corresponds to $\xi K$ and $\R e_2$ to $\nu K$. Then  
$n(k_1)=n(0)=e_1$ and $n(k_2)=n(\pi)=-e_1$. Define 
$n(\theta)=\cos(\theta)e_1+\sin(\theta)e_2$, $0\le\theta\le\pi$.   
\end{itemize}

We use $C$ to compute the linking number $\lk(K,K_n)$. 
We first prove the theorem under the assumption that all components
$K$ of $V$ are contractible:  

Let $W$ be a sufficiently small neighborhood of $\R V$ in $\C P^3$ and
let $U=\R P^3\cap W$, as in Remark ~\ref{rmknearfar} (we use notions
$C(U)$, $C'(U)$, $\Gamma(W)$, and $\Gamma'(W)$, as there).
First, we relate $\wra(\R V,n)$ to the sum of local writhes at double
points $x$ of $p_c(\R V)$ (the double points labeled (\rdp) in
\S ~\ref{diagapp}). Since the vector field $n$ is everywhere orthogonal
to $TC$, $C(U)\cap K_n=\emptyset$. 

To each double point $x$ of type (\rdp) there corresponds
exactly two points in $C'(U)\cap K$, which are the preimages of
$x$ under $p_c$. Denote these two points $a$ and $b$, respectively. We
use the notions $v$, $w$, $w'$, $S$, and $u$ as in \S
~\ref{diagapp}. Extend $u$ to a nonzero tangent vector 
field of $S$. Then, according to our orientation convention and since
$\pa C=2K$, $(u,w')$ is a positively oriented basis of
$TC'(U)$ at $a$ and the local intersection number equals the sign of the
orientation $(v,u,w')$, which is just the local writhe.  Interchanging
$a$ and $b$ in the above argument, we 
see that the local intersection number of $C'(U)$ and $K$ at
$b$ also equals the local writhe.  
Hence, the sum of local writhes at double points of type (\rdp) the
preimages of which lies on $K$, equals $\lk(K,K_n)$. It follows that
the sum of local writhes at double points of type (\rdp), the preimages
of which lies on the same component, equals $\wra(\R V,n)$.

Second, we relate $\sh(V,n)$ to the sum of local
writhes at double points of type (\sdp). Use $c$ to construct the shade
$\Gamma=\Gamma_c$. Since $n$ is orthogonal to $TC$, it follows that 
$i n$ is transverse to $T\Gamma_c$ along $\R V$. Hence,
$\Gamma(W)\cap\C V_n=\emptyset$. Then
$$
\Gamma\bullet\C V_n=
\Gamma'(W)\bullet\C V_n=
\Gamma'(W)\bullet\C V.
$$

To each double point $s\in \R H$ of $P_c(\C  V)$ (of type (\sdp)) there
corresponds 
exactly two complex conjugate points in $\Gamma'(W)\cap\C V$, which are
the preimages of 
$s$ under $P_c$. Denote these two points $a$ and $b$, respectively.
Note that $a$ and $b$ lie in 
$\C l\subset\Gamma$ where $l$ is the line through $s$ and $c$. 

Consider the local intersection number of $\C V$ and $\Gamma'(W)$ at
$a$. We use the notions  $\D l(a)$, $\D l(b)$, $\C V(a)$, and 
$\C V(b)$ as in \S ~\ref{diagapp}.    
Let $\alpha\subset\D l(a)$ be an arc connecting $a$ to $s$.
Push $\C V(a)$ along $\alpha$ to $s$, keeping it transverse to $\Gamma$,
and so that at the end of the push $\C V(a)$ agrees (locally around $s$)
with $P_c(\C V(a))$. From this construction it follows that 
the local intersection number of $\Gamma'(W)$ and $\C V(a)$ at $a$ equals 
the local intersection number of $\Gamma$, with orientation
coming from the part of $\Gamma$ containing $\D l(a)$, and 
$P_c(\C V(a))$ at $s$.  

The orientation of $T_s\Gamma$  coming from the part of
$\Gamma$ containing $\D l(a)$ (see \S
~\ref{shadeofRP} and ~\ref{shadeofRS}) is given by the orientation of
the frame 
$(v,iv,u,w)$, where $v$ is a real tangent vector along $\R l$
pointing in the direction of the orientation of $l$ induced from
$\D l(a)$, and $(u,w)$ is a frame in $T_sH$ such that $(u,w,v)$ is a
positive basis of $T_s\R P^3$. 

Let $(f,if)$ be a complex frame of $T_sP_c(\C V(a))\subset T_s\C H$
(neither $f$ nor $if$ lies in $T_s \R H$). Then the local
intersection number of $\Gamma$ and $P_c(\C V(a))$ is given by the
sign of the orientation of the frame 
$$
(v,i v,u,w,f,if),
$$  
which is positive if and only if the frame $(u,w,f,if)$ gives
the complex orientation of $\C H$.

On the other hand, the local writhe of $s$ is the sign of the
orientation of the frame $(v,u',w')$ of $T_s\R P^3$, where
$(u',w')$ is a frame of $T_s\R H$ 
such that $(u',w',f,if)$ gives the complex orientation of 
$\C H$. So, the local writhe is positive if and only if 
$(u,w,f,if)$ gives the complex orientation of $\C H$.

It follows that the sum of local writhes at double points of
type (\sdp) equals $\sh(V,n)$.

The above together with Theorem ~\ref{thmodd} proves the theorem for
real algebraic links without non-contractible components.

Now assume that $\R V$ has non-contractible components. The only
difference from the case considered above is that, for each
non-contractible component $K$, there will be one intersection point
of $C(U)$ and $K_n$ (a point of type (\rn)) and a corresponding
intersection point of $\Gamma(W)$ and $\C V_n$ (of type (\cn)). These
intersection points 
appear near the sub-arc $J$ of $K$. We compare their signs (cf. the
end of Remark ~\ref{rmksigns}).   

We use local coordinates in a
neighborhood of $J\subset K$ as in the proof of Lemma
~\ref{lmashadetw} (the case $k=1$):
$K$ is given by the equation $x_2=x_3=0$,  
$C\simeq A=a_1\cup a_2$, and $\Gamma\simeq\gamma_1\cup\gamma_2$. 

We may assume that $0\le x_1\le\pi$ corresponds to $J$, that the
orientation of $K$ is the one from $0$ to $\pi$, and that $n(0)=\partial_2$.
The intersection point of $A$ and $K_n$ (type (\rn)) is the point
$(\frac{\pi}{2},0,\delta)$ and the sign of the intersection is the sign of
the frame $(\partial_1,\partial_3,-\partial_2)$, which is $+1$.
The intersection point of $\Gamma$ and $\C V_n$ (type (\cn)) is
$(\frac{\pi}{2},0,i\delta)$ and the sign of the intersection is the sign of
the frame
$$
(\partial_1,\partial_2,\partial_3,i\partial_3,-i\partial_2,i\partial_1),
$$
which is $-1$. Hence, the contributions to $\EW$ from the two extra
intersection points  
for each non-contractible component cancel and the theorem follows
in the general case.  
\qed


\section{Generalizations}
In this section, it is shown that any  
$k$-dimensional subvariety of $(2k+1)$-\linebreak
dimensional complex projective
space or of the complexification of the real $(2k+1)$-sphere can be
equipped with additional structure 
which allows for shade numbers to be defined. Throughout this section, 
we call the complexification of the real $(2k+1)$-sphere the complex
$(2k+1)$-sphere and
we let $\R Y^{2k+1}$ denote $\R P^{2k+1}$ or $\R Q^{2k+1}$ and 
$\C Y^{2k+1}$ be the associated complex manifold

\subsection{Additional structure}\label{secgeneral}
Let $W$ be a projective $k$-dimensional variety in complex projective
$(2k+1)$-space or in the complex $(2k+1)$-sphere.
A vector field $n$ in $\R Y^{2k+1}$ along $\R W$ (i.e. a section of
the bundle $T\R Y^{2k+1}$ restricted to $\R W$) is {\em admissible}
if it has the following properties.
\begin{itemize}
\item $n$ is the restriction of a smooth vector field defined in some
neighborhood of $\R W$ in $\R Y^{2k+1}$.
\item For all smooth extensions $\nu$ of $in$ supported in 
sufficiently small neighborhoods of $\R W$ in $\C Y^{2k+1}$, the flow
$\Phi_t$ of $\nu$ satisfies $\Phi_t(\C W)\cap\R Y^{2k+1}=\emptyset$,
for all sufficiently small $t>0$.
\end{itemize}

If $(W,n)$ is a variety with an admissible vector field $n$ then its
shade number can be defined as follows. 
\begin{equation}\label{dfngeneral}
\sh(W,n)=\frac12([\Gamma]\bullet[\C W_n])\in\frac12\Z,
\end{equation}
where $[\C W_n]$ denotes the homology class in 
$H_{2k}(\C Y^{2k+1}-\R Y^{2k+1})$ of the cycle $\C W_n$ obtained by
shifting 
$\C W$ slightly by the flow of $\nu$. The shade number is
independent of the extension $\nu$ of $in$ as long as its support is
sufficiently small. 

A general position argument shows that any $k$-dimensional variety
has an admissible vector field: 
$\C W$ is a stratified set of real dimension $2k$ and standard
transversality arguments applied to a tubular neighborhood of 
$\R Y^{2k+1}$ in $\C Y^{2k+1}$ show that there exists a smooth purely
imaginary vector field $v$ along   
$\R Y^{2k+1}$ in $\C Y^{2k+1}$ supported in an arbitrarily small
neighborhood of $\R W$, such that if $\Phi_t$ is any
$1$-parameter deformation of $\R Y^{2k+1}$ with
$\frac{d}{dt}\Phi_t|_{t=0}=v$ then 
$\Phi_t(\R Y^{2k+1})\cap\C W=\emptyset$, for all sufficiently small
$t>0$.  
Let $n(x)=iv(x)$ for $x\in\R W$. Then $in=-v$ and $n$ is admissible.

The cases considered earlier in this paper, varieties with empty real
set (\S ~\ref{IRempty}) and armed real varieties without real
singularities (\S ~\ref{Ishadearm}) are special cases of the 
definition given in Equation ~\eqref{dfngeneral}. (In fact, they are the
generic cases for complex respectively real varieties.) In particular, the
case of armed real varieties shows that admissible vector fields are not
unique and that different choices may give different shade numbers. 

We now turn our attention to some particular non-generic cases which
are similar to the case of real varieties studied earlier in that
admissible vector fields are easy to find and their non-uniqueness
can be controlled. 


\subsection{Varieties with manifolds as real sets}\label{secRVmfd}
Let $W$ be a projective $k$-dimen-sional variety in complex projective
$(2k+1)$-space or in the complex $(2k+1)$-sphere. Assume that
$\R W$ is submanifold of $\R Y^{2k+1}$ and that the points in 
$\R W$ are smooth points of $W$. (The connected components of 
$\R W$ may have different dimensions, the maximal possible dimension
of a component is $k$.) Assume also that if $x\in\R W$ is any point in an
$r$-dimensional component $K$ of $\R W$ then 
$T_x\C W\cap T_x\R P^{2k+1}= T_x K$. (In other words, $\C W$ and 
$\R Y^{2k+1}$ intersect cleanly along the manifold $\R W$).  

Under these conditions, an admissible vector field along $\R W$ can
be constructed as follows. 

Let $K$ be an $r$-dimensional component of $\R W$. The normal 
bundle $N K$ of $K$ in $\R Y^{2k+1}$ is $(2k+1-r)$-dimensional. It
has a natural $(r+1)$-dimensional subbundle $N'K$: if $x\in K$
then $N'_xK=T_x\R Y^{2k+1}/\pr(T_x\C W)$, where 
$\pr\colon T_x\C Y^{2k+1}\to T_x\R Y^{2k+1}$ is the projection. For
dimensional reasons there exists non-zero sections of
$N'K$. If $n$ is a vector field along $\R W$ which projects to a
non-zero section of $N' K$ for each component $K$ of $\R W$ then
$n$ is admissible. 

Conversely, any admissible vector field gives rise to such
sections. Moreover, two admissible vector fields which give rise to 
homotopic sections in $N'K$ for every $K$ clearly give the
same shade number. It is straightforward to describe all possible
homotopy classes of non-zero sections in these bundles.

If $K=\{x\}$ is a $0$-dimensional component of $\R W$, then the 
bundle $N'K$ is $1$-dimensional. The tangent space $T_x\R Y^{2k+1}$
has a preferred orientation from a fixed orientation of $\R Y^{2k+1}$
and the image $\pr(T_x\C W)$ has an orientation induced by the
complex orientation of $T_x\C W$.
This allows us to define $n(x)$ in a canonical way:
let $n(x)$ be the {\em positive normal} of $\pr(T_x W)$ in
$T_x\R Y^{2k+1}$. It is easy
to check that each $0$-dimensional component of $\R W$ endowed with the
canonical vector field contributes $(-1)^k\frac12$ to $\sh(W,n)$. This
observation leads to the following analogue of Theorem ~\ref{thmrange} for
varieties $W$ in $\C P^{2k+1}$ with $\R W$ $0$-dimensional:

\begin{thm}\label{thmrange2}
For varieties $W$ as above of degree $d$ and with
$\R W=\{p_1,\dots,p_m\}$ (note that $m\le d^2$), the range of the
shade number consists of all half-integers between
$(-1)^k\frac12(m-d^2)$ and $(-1)^k\frac12d^2$ congruent to $\frac12d$
modulo $1$. \qed 
\end{thm}

\subsection{Real varieties with singularities} 
Let $V$ be a real projective $k$-dimensional variety in real projective
$(2k+1)$-space or in 
the real $(2k+1)$-sphere and let $n$ be an
admissible vector field along $\R V$. When $\dim(V)$ was odd, $V$ was
without real singularities, and $\R V$ orientable, we compared the
shade number to 
the wrapping number. One may study related issues in more general
situations there are however differences between the singular and the
non-singular cases: 

The shade number is defined once an admissible vector field
has been picked. To have a counterpart of the  wrapping number
defined, there are two obvious requirements which must be met: 
that all small shifts of $\R V$ along $n$ shifts $\R V$ off
itself, and that $\R V$ is an orientable $k$-cycle in $\R P^{2k+1}$.
Even with these two conditions met, it is not clear how to define the
counterpart of the wrapping number of the $k$-cycle. One reason is
that an orientable connected $k$-cycle, in contrast to an orientable
connected manifold, may have more than two orientations. To get a
reasonable counterpart of the wrapping number one 
must fix a semi-orientation (an orientation up to sign) on the
$k$-cycle $\R V$. 

We next study properties of shade numbers for real varieties with the
simplest singularities: varieties with double points as that of $V_0$ in
Theorem ~\ref{thmjump} (this refers to the double points, $\R V_t$
need not have global properties as there, i.e. it need 
neither be connected nor orientable).  

In this case, admissible vector fields are easy to find: pick any
normal vector field along $\R V$ which is 
transverse to the (real) $2k$-dimensional intersection of the Zariski
tangent space of $V$ (i.e. of $V$ after base extension) and 
$T_x\R Y^{2k+1}$ at the double points $x$.  

In the odd-dimensional case ($\dim(V)$ is odd), the shade number 
of a variety with one double point is the mean value of the
shade numbers of its (two) resolutions, with the induced admissible
vector field, (if the double point is real-real then the shade numbers of
the resolutions are the same, if it is complex-complex-conjugate then
the shade numbers of the resolutions differ by $\pm2$). 

In the even-dimensional case, the shade number of a
variety with one complex-complex-conjugate double point differs from
the common value of the shade numbers of its resolutions by 
$\pm 1$, depending on the direction of the shift at the double point,
and the shade number of a variety with a real-real double point agree
with the common value of the shade numbers of its resolutions.


\section{Examples}
In this section, several examples are presented. 


\subsection{Two families of real algebraic knots}
Let $[x_0,x_1,x_3,x_4]$ be projective coordinates on real projective
$3$-space and let $a\in\R$. Let $\epsilon=\pm1$ and let
$K_a(\epsilon)$ be the variety 
defined by the equations 
\begin{align*}
& ax_0x_2-x_1x_3=0,\\
& a^2x_1x_0-\epsilon a^2x_0^2-x_3^2=0.\\
& x_3x_2^2+\epsilon x_1^2x_3-ax_1^2x_2=0,\\
& x_2^2x_0-x_1^3+\epsilon x_1^2x_0=0.
\end{align*}
It is easy to check that for $a\ne 0$, $K_a(\epsilon)$ are smooth
curves. The curve $K_0(\epsilon)$ has one 
double point. For $\epsilon=-1$ it is a real-real double point and for
$\epsilon=+1$ it is a complex-complex-conjugate double point. In fact,
$K_a(\epsilon)$ is the rational curve given, in projective
coordinates, by   
$$
[s,t]\mapsto\left[s^3,st^2+\epsilon s^3, t^3 +\epsilon s^2t,
ats^2\right]. 
$$

It follows from Theorem ~\ref{thmodd} (ii) that $\EW(K_a(\epsilon))$ changes
by $\pm 2$ at $a=0$. Hence, $\EW(K_a(\epsilon))\ne 
\EW(K_{-a}(\epsilon))$, $a\ne 0$ and 
thus $K_a(\epsilon)$ is not weak rigid isotopic to
$K_{-a}(\epsilon)$. However, it is clear  
that $\R K_{a}(\epsilon)$ and $\R K_{-a}(\epsilon)$ are topologically 
isotopic. Also the projective isomorphism
$[x_0,x_1,x_2,x_3]\mapsto[x_0,x_1,-x_2,x_3]$ takes  
$K_{a}(\epsilon)$ to $K_{-a}(\epsilon)$.      


\subsection{An armed real projective plane}\label{RP2}
Let $[x_0,x_1,x_2]$ and $[y_0,\dots,y_5]$ be projective coordinates 
on real projective $2$- and $5$-space, respectively.  
Consider the map $\phi_0\colon\R P^2\to\R P^5$ 
\begin{equation*}
\phi_0([x_0,x_1,x_2]) = [x_0,x_1,x_2,0,0,0]
\end{equation*}
Then $\phi_0$ gives a parameterization of the variety $V$ defined by
the equations $y_3=y_4=y_5=0$. Let $\phi_t$ be the $1$-parameter
variation of $\phi_0$ given by  
\begin{equation*}
\phi_t([x_0,x_1,x_2])= [x_0,x_1,x_2,tx_0,tx_1,tx_2],
\end{equation*}
and let $n=\frac{d}{dt}|_{t=0}\phi_t$. Then $n$ is a normal vector
field along $\R V^2$ in $\R P^5$. 

\begin{clm}
If the armed variety $(V,n)$ is as above 
then $\sh(V,n)=-\frac12$. 
\end{clm}

\begin{pf}
Consider $\Gamma_p$, where $p=[0,0,0,0,0,1]$. It is
straightforward to check that  
$\C V_n\cap\Gamma_p=\{[0,0,1,0,0,i\delta]\}$, where $\delta>0$ is
small. In coordinates 
$$
(x_1+i y_1,\dots,x_5+iy_5)\mapsto[1,x_1+i y_1,\dots,x_5+iy_5],
$$ 
the sign of the intersection point equals the sign of the orientation
of the frame
$$
(\pa_1,\dots,\pa_5,i\pa_5,\pa_1+i\pa_3,i\pa_1,\pa_2+i\pa_4,i\pa_2),
$$  
which is $-1$.
\end{pf}


\subsection{An unknot}\label{Oknot}
Let $[x_0,\dots,x_4]$ be homogeneous coordinates on real projective
$4$-space. Let $Q^3$ be as in \S ~\ref{shadeofRS}. If
$O$ is the intersection of $Q^3$ and the $2$-plane given by 
$x_3=x_4=0$ then $\R O\subset\R Q^3$ is a representative
of the unknot. 

\begin{clm}
If $O$ is as described above then $\EW(O)=0$.
\end{clm}

\begin{pf}
Let $[s,t]$ be projective coordinates on the real projective line. The
variety $O$ admits the rational parameterization 
\begin{equation*}
[s,t]\mapsto\left[s^2+t^2,2st,s^2-t^2,0,0\right].
\end{equation*}

Let $p=[0,0,0,0,1]\in\R P^3$ (as in \S ~\ref{shadeofRS}
we think of $\R P^3\subset\R P^4$ given by the equation $x_0=0$). 
Consider $\Gamma=\Gamma_p$ and $C$, where 
$C\subset\R Q^3$, $\partial C=2O$, is the chain constructed
using $\tilde p_0=[1,0,0,0,1]$, see \S ~\ref{generic}. 

Let $W$ be a small neighborhood of $\R O$ in $\C Q^3$ and let
$U=W\cap\R Q^3$, as in Remark ~\ref{rmknearfar} (and using
the same notation as there), $\EW(O)$ equals half of the sum of the 
intersection numbers $\Gamma'(W)\bullet\C O$ and $C'(U)\bullet O$ (see 
the proof of Theorem ~\ref{thmodd} and Remark ~\ref{rmksigns}).  

It is straightforward to check that $\R O\cap C'(U)=\emptyset$.
A point 
\begin{equation*}
[x_0,x_1,x_2,0,0]=\left[s^2+t^2,2st,s^2-t^2,0,0\right]\in \C O
\end{equation*}
(here we think of $[s,t]$ as homogeneous coordinates on $\C P^1$)
lies in $\Gamma'(W)$ if and only if it is not a real point and all the
products $x_ix_j^{\ast}$, $i\ne j$ are real numbers ($x^\ast$ denotes the
complex conjugate of $x$). Since $[0,1]\in\C P^1$ maps to a real point
we can restrict attention to the coordinate chart 
$t\in \C\mapsto [1,t]\in\C P^1$. Then $x_1x_2^{\ast}\in \R$ only if
$\rho 2t=1-t^2$ for some real number $\rho$. However, this equation has 
only real solutions for every $\rho\in\R$. Hence 
also $\Gamma'(W)\cap\C O=\emptyset$, and we conclude $\EW(O)=0$. 
\end{pf}


\subsection{A trefoil knot}\label{32knot}
Let $[x_0,\dots,x_4]$ be coordinates on real projective $4$-space and
let $(z,w)$ be coordinates on $\C^2$.  The trefoil knot is the link of
the singularity $z^2=w^3$. We think of $\C^2\approx\R^4$ as the affine
part of $\R P^4$ given by the condition $x_0\ne 0$. For simpler formulas
below we violate our conventions slightly and represent $\R Q^3$ as
the quadric $-2x_0^2+x_1^2+\dots+x_4^2=0$. The intersection of
the subset $z^2=w^3$ of $\R^4\subset\R P^4$ and $Q^3$ is 
a variety $K$ such that $\R K\subset\R Q^3$ is a
representative of the trefoil knot. 

\begin{clm}
If $K$ is as above then $\EW(K)=4$. 
\end{clm}

\begin{pf}
The knot $\R K$ can be parameterized by $\theta\in[0,2\pi]\subset\R$ as 
\begin{equation*}
\theta\mapsto\left[1,\cos\,3\theta,\sin\,
3\theta,\cos\,2\theta,\sin\,2\theta\right].  
\end{equation*}
Using the rational parameterization
$[s,t]\mapsto\left(\frac{2st}{s^2+t^2},\frac{s^2-t^2}{s^2+t^2}\right)$
of $\cle^1$ 
we get the rational parameterization
$\phi([s,t])=[x_0([s,t]),\dots,x_4([s,t])]$  of $K$, where 
\begin{align*}
x_0([s,t]) &= (s^2+t^2)^3,\\
x_1([s,t]) &= -2st\left(3(s^2+t^2)^2-16s^2t^2\right),\\
x_2([s,t]) &= -(s^2-t^2)\left((s^2+t^2)^2-16s^2t^2\right),\\
x_3([s,t]) &= -(s^2+t^2)\left((s^2+t^2)^2-8s^2t^2\right),\\
x_4([s,t]) &= 4st(s^2+t^2)(s^2-t^2).
\end{align*}   
(It is easy to check that the same formulas over complex numbers give a
holomorphic immersion of $\phi\colon\C P^1\to\C Q^3$ and thus $\phi$
parameterizes $K$.)

We use the points $\tilde p_0=[1,0,\sqrt{2},0,0]$ and 
$\tilde p_1=[1,0,-\sqrt{2},0,0]$ to construct 
the shade $\Gamma\subset\C Q^3$, and
the chain $C\subset\R Q^3$, $\partial C=2K$, see \S
~\ref{generic}. The sphere $\R Q^3$ has the standard 
orientation coming from $\C^2\subset\R P^4$. As in \S
~\ref{Oknot}, we must calculate $C'(U)\bullet \R K$ and 
$\Gamma'(W)\bullet\C K$, where $W$ is a small neighborhood of $\R K$
in $\C Q^3$ and $U=\R Q^3\cap W$. 

We start with $C'(U)\bullet \R K$. Using the parameterization of $\R K$ by
$\theta\in[0,2\pi]\subset\R$ as given above, there are nine point
pairs in $C'(U)\cap \R K$. These nine pairs $(\alpha,\beta)$ are of
two types:
\begin{equation*}
(\alpha,\beta)=\begin{cases}
\text{type 1: }(\frac{\pi}{4}+n\frac{\pi}{3},\frac{3\pi}{4}+n\frac{\pi}{3}) &
n=0,1,2,3,4,5,\\
\text{type 2: }
(\frac{\pi}{6}+n\frac{\pi}{3},\frac{7\pi}{6}+n\frac{\pi}{3}) &
n=0,1,2.
               \end{cases}
\end{equation*}   
Each point pair contributes $\pm 2$ to $C'(U)\bullet \R K$. Using
stereographic projection $s\colon\R Q^3-\{\tilde p_0\}\to\R^3$,
$s(\tilde p_1)=0\in\R^3$ as in \S ~\ref{pfthmodd}, one finds that
a point pair of type 1 contributes $(-1)^n2$ to $C'(U)\bullet \R K$ and that
the point pair of type 2 contributes $+2$. Hence,
$\frac{1}{2}(C'(U)\bullet \R K)=3$.

We proceed to calculate $\Gamma'(W)\bullet \C K$. Note that $\phi([s,t])$
lies in $\Gamma$ if and only if $x_i([s,t])x_j^\ast([s,t])\in \R$, 
$i,j\in\{1,3,4\}$ ($x^\ast$ denotes the complex conjugate of $x$). We
are looking for points in $\Gamma'(W)\cap\C K$, and 
so we are interested in non-real points $[s,t]\in\C P^1$ for which
this condition holds. Since $[1,0]$ and $[0,1]$ are real points in 
$\C P^1$ we may work in the coordinate chart
$t\in\C-\{0\}\mapsto[1,t]\in\C P^1-\{[1,0],[0,1]\}$.  
The condition $x_1([1,t])x^\ast_4([1,t])\in\R$ is then
equivalent to 
\begin{equation*}
(3(1+t^2)^2-16 t^2)(1-(t^\ast)^4)\in\R,
\end{equation*}
which holds if and only if $t=\pm 1,\pm i$ or 
\begin{equation*}
(3(1+t^2)^2-16 t^2)=\rho(1-t^4), \text{ for $\rho\in\R$}.
\end{equation*}
This equation has solutions $t^2=\frac{5\pm\sqrt{16+\rho^2}}{3+\rho}$
if $\rho\ne -3$ and $t=\pm\frac{3}{\sqrt{10}}$ if $\rho={-3}$. We
conclude that $x_1([1,t])x^\ast_4([1,t])\in\R$ only if $t\in\R$ or
$t\in i\R$. The latter case is the one of interest. For $t=ia$,
$a\in\R-\{0\}$ the condition $x_4([1,t])x_3^\ast([1,t])\in \R$ reads
\begin{equation*}
4ia(1+a^2)(1-a^2)^2((1-a^2)^2+8a^2)\in\R.
\end{equation*}
This condition is met for $a=\pm 1,\pm i$ or $a^2=-3\pm\sqrt{8}$ and
hence the only points in $\R-\{0\}$ which meet the condition is 
$a=\pm 1$. We conclude that the only points in $\Gamma'(W)\cap \C K$ are
$\phi([1,\pm i])=[0,\pm i, 1,0,0]$. We must calculate the intersection
number. Around $[0,-i,1,0,0]$, the shade $\Gamma$ can be parameterized
by $(z,h,k)\in U\times I\times I$, where $U\subset\C$ is a small disk
around $i$ and $I\subset\R$ denotes a small open interval around $0$,
as follows 
$$
(z,h,k)\mapsto
\left[\frac{z^2+1}{\sqrt{2}(z^2-1)},\frac{2z\sqrt{1-h^2-k^2}}{z^2-1},
1,\frac{2hz}{z^2-1},\frac{2kz}{z^2-1}\right],
$$ 
and the orientation of $\Gamma$ is the one induced by the complex
orientation of $U$ followed by the orientation of the frame
$(\partial_h,\partial_k)$ on $I\times I$. A straightforward
calculation shows that the intersection number of $\Gamma'(W)$ and 
$\C K$ at $[0,1,-i,0,0]$ is $+1$. By the proof of Theorem ~\ref{thmodd}
the intersection number at $[0,1,i,0,0]$ is the same. Thus
$\frac{1}{2}(\Gamma'(W)\bullet \C K)=1$. Adding the terms,
$\EW(K)=3+1=4$.  
\end{pf}


\end{document}